\documentclass[12pt]{article}
\usepackage[latin1]{inputenc}
\usepackage{amsmath,amssymb}
\usepackage{graphics, graphicx,psfrag}
\usepackage{amsthm}
\usepackage{epsfig}
\usepackage{amscd}
\usepackage{verbatim} 
\usepackage[english,frenchb]{babel}

\newcommand{\R}{R}
\newcommand{\Q}{Q}  

\newcommand{\G}{\Gamma}
\renewcommand{\S}{S}

\newcommand{\I}{\mathrm{I}}
\newcommand{\IIi}{\mathrm{II_1}}
\newcommand{\II}{\mathrm{II}}
\newcommand{\III}{\mathrm{III}}

\newcommand{\CI}{{\mathrm{\bf C}}}
\newcommand{\RI}{{\mathrm{\bf R}}}
\newcommand{\QI}{{\mathrm{\bf Q}}}
\newcommand{\ZI}{{\mathrm{\bf Z}}}

\newcommand{\SI}{{\boldsymbol S}}

\newcommand{\SL}{\mathrm {SL}}

\newcommand{\Sp}{\mathrm{Sp}}

\newcommand{\norm}[1]{||#1||}
\newcommand{\abs}[1]{|#1|}

\newcommand{\1}{{\bf 1}}
\newcommand{\graph}{\mathrm{graph}}

\newcommand{\impl}{\Longrightarrow}
\newcommand{\ssi}{\Longleftrightarrow}

\newcommand{\surj}{\twoheadrightarrow}

\newcommand{\del}{\partial}
\renewcommand{\limsup}{\overline{\mathrm{lim}}\ }

\newcommand{\id}{\mathrm{Id}}

\renewcommand{\Re}{\mathop{\mathrm{Re}}}

\newcommand{\h}{{\mathfrak{h}}} 

\newtheorem*{dfn}{\indent\bf D\'efinition}
\newtheorem{thm}{\indent\bf Th\'eor\`eme}
\newtheorem{lem}[thm]{\indent\bf Lemme}

\newtheorem{prop}[thm]{\indent\bf Proposition}
\newtheorem{cor}[thm]{\indent\bf Corollaire}

\renewcommand{\qed}{_\blacksquare}
\newenvironment{dem}{\noindent D\'emonstration. }{{\hfill {$\qed$}}\vspace{.25cm} }

\makeatletter

\renewcommand\thepart{\@Roman\c@part.}
\renewcommand\thesection{\@arabic\c@section.}
\renewcommand\thesubsection{\@arabic\c@subsection.}

\renewcommand\section{\@startsection {section}{2}{\z@}%
    {0ex \@plus -1ex  \@minus -.2ex}%
    {2.3ex \@plus .2ex}%
    {\reset@font \center \sc }}
\renewcommand\subsection{\@startsection {subsection}{3}{\parindent}%
    {.1ex \@plus .5ex \@minus .2ex}%
    {-.8em \@plus .2ex}%
    { \reset@font\small \bfseries}}

\makeatother

\usepackage{diagrams}
\usepackage{array}

\begin{document}

\centerline{\bf{\uppercase{Sur les espaces mesur\'es singuliers}}}

\medskip

\centerline{{\bf II - \'Etude spectrale}\footnote{
The author was partially supported by a JSPS Fellowship for European Researchers.}}
\vspace{1cm}
\centerline{\small {Mika\"el \sc Pichot}}
\vspace{.4cm}
\centerline{\small{Novembre 2004}}
\vspace{.2cm}
\centerline{---}
\vspace{.3cm}

\selectlanguage{english}
\begin{abstract}
We are mainly interested here in  Kazhdan's property T for measured equivalence relations.   Among our main results are characterizations of strong ergodicity and Kazhdan's property in terms of the spectra of  diffusion operators, associated to random walks and hilbertian representations of the underlying equivalence relation. The analog spectral characterization of property T for countable groups was proved recently by Gromov  \cite{Gromov03} (and Ghys \cite{Ghys03}). Our proof put together the tools developed in the group case and further crucial technical steps from the study of amenable equivalence relations in \cite{CFW81}. As an application we show how \.Zuk's ``$\lambda_1 >1/2$" criterion for property T can be adapted to measured equivalence relations.
\end{abstract}
\selectlanguage{french}

\tableofcontents\bigskip\bigskip\bigskip

\section{Introduction}

Cet article succ\`ede \`a \cite{Pichot04_I} et concerne principalement l'\'etude de la propri\'et\'e T de Kazhdan pour les espaces singuliers. 

\bigskip

Nous avons divis\'e cette introduction en trois parties, une pr\'esentation des r\'esultats spectraux \'evoqu\'es dans le titre, suivie de quelques motivations, historiques et techniques, concernant  ces r\'esultats et  la propri\'et\'e T de Kazhdan ; nous terminons par exemple d'application.

\bigskip
\bigskip
\centerline{---}
\bigskip

\'Etant donn\'e un espace bor\'elien standard  $X$, on construit une relation d'\'equivalence bor\'elienne \`a classes d\'enombrables en effectuant une marche al\'eatoire sur $X$ de la fa\c con suivante.  Fixons pour tout $x\in X$ une loi de probabilit\'e $\nu_x$ sur $X$, \`a support fini ou d\'enombrable. Fixons un point  $x_0\in X$ et lan\c cons une marche al\'eatoire sur $X$ d\'ebutant en $x_0$ et de loi $\nu$. Ainsi on effectue un premier pas en $x_1$ selon la loi $\nu_{x_0}$, un deuxi\`eme  pas en  $x_2$, selon la loi $\nu_{x_1}$, et on it\`ere. L'ensemble des trajectoires possibles de cette marche d\'etermine une orbite d\'enombrable $[x_0] \subset X$ (l'ensemble des points  atteints \`a partir de $x_0$ avec probabilit\'e non nulle), et lorsque le point de d\'epart $x_0$ varie dans $X$, la partition en orbites de $X$ ainsi d\'etermin\'ee est une relation d'\'equivalence $\R$ sur $X$ (o\`u l'on suppose que $\nu_x(y)>0$ si et seulement si $\nu_y(x)>0$). Cette relation d'\'equivalence est bor\'elienne si le support  $\mathrm{Supp}(\nu) =\{  (x, y),\  \nu_x(y)>0\}$ de $\nu$ est une partie bor\'elienne de $X\times X$. On notera $\nu(x\to y)=\nu_x(y)$.

\bigskip

Fixons alors un espace mesur\'e singulier $\Q$. Rappelons que la structure mesur\'ee sur $\Q$ est d\'etermin\'ee par la donn\'ee d'un espace de probabilit\'e $(X,\mu)$ et d'une relation d'\'equivalence mesur\'ee $\R$ \`a classes d\'enombrables sur $X$ telle que $X/\R = \Q$. (On dit alors que $\R$ est une relation d'\'equivalence d\'esingularisante relative \`a $\Q$.) On appelle  \textit{marche al\'eatoire $\Q$-p\'eriodique}, et on note $\tilde \nu$, la donn\'ee d'une marche al\'eatoire $\nu$ sur un espace de probabilit\'e $(X,\mu)$ (au sens ci-dessus), d\'efinissant une relation d'equivalence bor\'elienne $\R$ pour laquelle $\mu$ est une mesure \textit{quasi-invariante}, et telle que $X/\R=\Q$. La marche al\'eatoire $\Q$-p\'eriodique $\tilde \nu$ a lieu sur l'ensemble d\'enombrable $\Q$-p\'eriodique $\R$. (cf. \cite{Pichot04_I} pour la terminologie utilis\'ee)\\

Soit $\tilde \nu$ une marche al\'eatoire $\Q$-p\'eriodique. Certaines propri\'et\'es  de l'espace singulier $\Q$ se refl\`etent dans le spectre d'op\'erateurs de diffusion naturellement associ\'es \`a $\tilde \nu$ de la fa\c con suivante.  La relation d'\'equivalence $\R$, partition de $X$ en les orbites de $\nu$, est en particulier un groupo\"\i de mesur\'e, et admet \`a ce titre des repr\'esentations hilbertiennes. Soit $\pi$ une telle repr\'esentation sur un champ mesurable $H$ d'espaces de Hilbert de base $X$. On associe \`a $\nu$ et $\pi$ un \textit{op\'erateur de diffusion} $D_{\nu,\pi}$ agissant sur l'espace $L^2(X,\mu,H)$ des sections de carr\'e int\'egrable de $H$, d\'efini par 
\[
(D_{\nu,\pi}\xi)_x = \sum_{y\sim x} \nu(x\to y)\pi(x,y)\xi_y
\]
o\`u $\xi\in L^2(X,\mu,H)$. Lorsque $\nu$ est une marche al\'eatoire \textit{sym\'etrique} relativement \`a $\mu$ (cf. \textsection \ref{dem}), cet op\'erateur est hermitien born\'e et son spectre est une partie compacte de $[-1,1]$. On s'int\'eresse  alors \`a la propri\'et\'e suivante.

\bigskip

On dit que $D_{\nu,\pi}$ a \textit{un trou spectral} (au voisinage de 1) si sa plus grande valeur spectrale distincte de 1 est strictement inf\'erieure \`a 1.

\bigskip

Soit $\tilde \nu$ une marche al\'eatoire $\Q$-p\'eriodique.  On appelle \textit{diffusions hilbertiennes} associ\'ees \`a $\tilde \nu$ la famille des op\'erateurs $D_{\nu,\pi}$, lorsque $\pi$ parcourt les repr\'esentations hilbertiennes de la relation d'\'equivalence mesur\'ee $\R$ associ\'ee \`a $\tilde \nu$. On appelle \textit{diffusion simple} associ\'ee \`a $\tilde \nu$ la diffusion hilbertienne, agissant sur $L^2(X)$,  associ\'ee \`a la repr\'esentation triviale de $\R$.\\

Apr\`es avoir men\'e une analyse d\'etaill\'ee de la propri\'et\'e T de Kazhdan pour les espaces singuliers, nous obtenons dans cet article les caract\'erisations spectrales suivantes. 

\bigskip

\begin{thm}\label{ef-intr}
Un espace singulier ergodique de type fini $\Q$, s'il est de type $\II$, poss\`ede un quotient moyennable si et seulement si pour toute marche al\'eatoire $\Q$-p\'eriodique sym\'etrique born\'ee, la diffusion hilbertienne simple associ\'ee  n'a pas de trou spectral au voisinage de 1.
\end{thm} 

\medskip

\begin{thm}\label{T-intr}
Un espace singulier ergodique $\Q$ poss\`ede la propri\'et\'e T de Kazhdan si et seulement s'il existe une marche al\'eatoire $\Q$-p\'eriodique sym\'etrique dont les diffusions hilbertiennes pr\'esentent un trou spectral au voisinage de 1.
\end{thm}

\bigskip
\bigskip
\centerline{---}
\bigskip

1. L'\'etude des marches al\'eatoires en milieu g\'eom\'etrique (alg\'ebrique) re\c coit une attention consid\'erable depuis une cinquantaine d'ann\'ees. Le premier r\'esultat obtenu dans ce domaine est probablement le th\'eor\`eme de P\'olya (1921) sur la r\'ecurrence/transcience des marches al\'eatoires dans $\ZI^n$ (cf. \cite{Harpe00}).

\bigskip

Le comportement spectral proprement dit des marches al\'eatoires, notamment la sensibilit\'e du spectre aux structures g\'eom\'etriques (alg\'ebriques) sous-jacentes, a \'et\'e \'etudi\'e pour la premi\`ere fois en 1959 par Kesten. Il d\'emontrait alors le th\'eor\`eme suivant.

\bigskip

\begin{thm}[Kesten \cite{Kesten59_RWG,Kesten59_BM}] Soit $\G$ un groupe de type fini. On fixe une marche al\'eatoire invariante sym\'etrique $\tilde \nu$ sur $\G$, dont le support est un syst\`eme g\'en\'erateur sym\'etrique fini de $\G$, et on note $D_{\tilde\nu}$ l'op\'erateur de convolution associ\'e \`a $\tilde \nu$, agissant sur $\ell^2(\G)$. Alors  $D_{\tilde \nu}$ a un trou spectral si et seulement si $\G$ est non moyennable.
\end{thm}

\medskip

Ici l'op\'erateur $D_{\tilde \nu}=D_{\tilde \nu,\lambda}^\G$ est associ\'e \`a la marche al\'eatoire $\tilde \nu$ et \`a la repr\'esentation r\'eguli\`ere gauche $\lambda$ de $\G$  sur $\ell^2(\G)$. 

\bigskip

Plus r\'ecemment, il a \'et\'e observ\'e  par Gromov \cite{Gromov03,Ghys03} que l'on dispose  d'une caract\'erisation spectrale, analogue \`a celle de Kesten, de la propri\'et\'e T de Kazhdan (il s'agit d'une infime partie du contenu de \cite{Gromov03}). Un groupe de type fini $\G$ muni d'une marche al\'eatoire invariante sym\'etrique $\tilde \nu$ dont le support est un syst\`eme g\'en\'erateur sym\'etrique fini de $\G$, a la propri\'et\'e T de Kazhdan si et seulement si les diffusions $D_{\tilde \nu,\pi}^\G$ associ\'ees aux repr\'esentations unitaires $\pi$ de $\G$ ont un trou spectral. Cette caract\'erisation, bien qu'\'el\'ementaire,  est d'importance fondamentale. Elle permet par exemple de donner une  preuve conceptuelle du \og crit\`ere $\lambda_1 > 1/2$\fg\ pour la propri\'et\'e T (cf.  \cite{Ghys03,Ollivier03}). 

\bigskip
\bigskip

2. Une caract\'erisation spectrale analogue de la propri\'et\'e T de Kazhdan pour les relations d'\'equivalence mesur\'ees (et les espaces singuliers) est  plus difficile \`a obtenir ; nous ferons notamment un usage important  de techniques d\'evelopp\'ees par Connes-Feldman-Weiss \cite{CFW81} dans le cadre moyennable (un expos\'e r\'ecent et d\'etaill\'e des diff\'erents concepts de moyennabilit\'e fait l'objet de \cite{AnanRen01}). De fa\c con g\'en\'erale, nous devons d\'eterminer dans quelle mesure les propri\'et\'es spectrales des observables naturelles $D_{\nu,\pi}$  associ\'ees \`a $\Q$ sont sensibles au passage d'une notion de quasi-p\'eriodicit\'e $\Q$ triviale (i.e. d'une notion de p\'eriodicit\'e) \`a une notion de quasi-p\'eriodicit\'e non triviale --- la r\'eponse pouvant \textit{a priori} d\'ependre de la notion de quasi-p\'eriodicit\'e  choisie. Rappelons que, poursuivant le point de vue de \cite{Pichot04_I}, nous nous int\'eressons aux concepts de p\'eriodicit\'e (type $\I$) et quasi-p\'eriodicit\'e (type $\II$ et $\III$) qui sont d\'ecrits par les diagrammes suivants.

\begin{diagram}[nohug,notextflow]
\tilde Y &             &               &&&    & \tilde \Sigma      &      &   \\
\dTo     &\rdTo(2,1) & Y=\tilde Y/\G  &&&  & \dTo  & \rdTo(2,1) &  \Sigma=\tilde\Sigma/\R  &\\
        &    \ldTo_p     &            &&&  &  & \ldTo_p    &                         &\\
\Q=\{\star\} &         &              &&& & \Q &             &             \\
\mathrm{type} \ \I &         &        &&& &\quad \mathrm{type} \ \II\ \mathrm{et}\ \III  & & \\
\end{diagram}

\noindent Sur le diagramme de gauche, $\tilde Y$ est un complexe simplicial muni d'une action libre cocompacte d'un groupe d\'enombrable $\G$, et l'espace quotient $Y$ est un complexe simplicial compact. L'application $p$ est triviale et l'espace singulier $\Q$ est r\'eduit \`a  un point. Sur le diagramme de droite, $\tilde \Sigma$ est un complexe simplicial $\Q$-p\'eriodique. L'ensemble des sommets de $\tilde \Sigma$ est une relation d'\'equivalence mesur\'ee $\R$ agissant sur $\tilde \Sigma$, et l'espace quotient $X$ de $\tilde \Sigma$ par $\R$ est un espace bor\'elien standard muni d'une structure de lamination par complexes simpliciaux (et d'une classe de mesure d\'etermin\'ee par la mesure invariante sur $\Q$). L'application $p : X\to \Q$ est une d\'esingularisation simpliciale de $\Q$. L'espace singulier $\Q$ a le cardinal du continu. 

\bigskip
\bigskip

3. D\'ecrivons  plus pr\'ecis\'ement la situation. Le support d'une marche al\'eatoire $\Q$-p\'eriodique $\tilde \nu$ est un graphe $Q$-p\'eriodique $\tilde \Sigma$ et la marche al\'eatoire $\nu$ associ\'ee \`a $\tilde \nu$ a lieu sur la d\'esingularisation (discr\`ete) de $\Q$ constitu\'ee des sommets $X=\Sigma^{(0)}\subset \Sigma$ de la lamination $\Sigma$. Dans le cas des groupes, la pr\'esence de l'action cocompacte de $\G$  a un caract\`ere uniformisant et ram\`ene essentiellement les possibilit\'es de fluctuations de donn\'ees p\'erio\-diques (invariantes) \`a l'ensemble fini $X=Y/\G$. Ainsi les donn\'ees relatives \`a $\tilde \nu$ (notamment les probabilit\'es de transition $\tilde \nu(\gamma \to \gamma')$) sont en nombre fini. On peut alors facilement, par exemple, d\'eduire d'informations ponctuelles des contr\^oles uniformes sur ces donn\'ees.  Dans le cas quasi-p\'eriodique ce n'est plus le cas --- bien qu'on fasse toujours une hypoth\`ese de type \og cocompacit\'e\fg\ en requ\'erant que $\tilde \nu$ soit sym\'etrique born\'ee --- et il faudra  d\'eterminer si l'on peut s'affranchir de ces hypoth\`eses de cofinitude sans entra\^\i ner de modifications fondamentales sur les observations spectrales effectu\'ees. 

\bigskip

Nous verrons notamment que l'une des difficult\'es qui se posent provient du fait que, dans le cas quasi-p\'eriodique,  \textit{des comportements non triviaux peuvent appara\^\i tre sur des parties (quasi-p\'eriodiques) arbitrairement petites de l'espace} ; les contr\^oles  uniformes ont lieu sur une majorit\'e de l'espace seulement. Or ces \og fluctuations infinit\'esimales\fg\  peuvent avoir une influence sur  les propri\'et\'es spectrales des op\'erateurs $D_{\nu,\pi}$ (rappelons que les op\'erateurs $D_{\nu,\pi}$ agissent sur l'espace de Hilbert des \textit{sections de carr\'e int\'egrable} $L^2(X,H)$ associ\'ees au champ d'espaces de Hilbert $H$, o\`u la boule unit\'e de $L^2(X,H)$ contient des sections \`a support arbitrairement petit). La situation est r\'eminiscente de la construction des nombres de Betti $L^2$ pour les relations d'\'equivalence \cite{Gaboriau02}, o\`u l'auteur introduit des param\`etres de cut-off sur les homotopies pour en faire des op\'erateurs born\'es.

\bigskip
\bigskip

4. Consid\'erons d'abord le cas le plus simple de la repr\'esentation triviale (th\'eor\`eme \ref{ef-intr}).  Un exemple significatif pour lequel les fluctuations \'evoqu\'ees au paragraphe pr\'ec\'edent r\'ev\`elent effectivement une nature spectrale est donn\'e par l'existence de complexes simpliciaux quasi-p\'eriodiques contenant ou non  des suites de F\o lner \'evanescentes (cf. \cite{Pichot04_I}). Plus pr\'ecis\'ement, la pr\'esence de suites de F\o lner \'evanescentes dans le support de $\nu$ est d\'etect\'e par l'absence de trou spectral pour la diffusion associ\'ee \`a $\nu$ et \`a la repr\'esentation triviale (que nous avons appel\'ee \og diffusion simple\fg). Cette observation, qui repose essentiellement sur des techniques de Connes-Feldmann-Weiss \cite{CFW81} et Schmidt \cite{Schmidt81}, conduit au  th\'eor\`eme \ref{ef-intr} (cf. \textsection \ref{dem}).

\bigskip
\bigskip

5. Passons maintenant \`a la propri\'et\'e T de Kazhdan. Rappelons que, par d\'efinition, une relation d'\'equivalence mesur\'ee sur un espace de probabilit\'e $(X,\mu)$  poss\`ede la propri\'et\'e T si la condition suivante est v\'erifi\'ee.

\bigskip

\begin{tabular}{m{.7cm}m{12.5cm}}
(T) & Toute repr\'esentation hilbertienne poss\'edant une suite $\xi^n$ de champs de vecteurs presque invariante, unitaire au sens o\`u $\norm{\xi^n_x}_x = 1$ pour presque tout $x\in X$, poss\`ede un champ de vecteurs invariant \`a support total.
\end{tabular}

\bigskip

Nous rappelerons en d\'etail la terminologie utilis\'ee ult\'erieurement ; nous voulons seulement ici attirer l'attention  sur la signification de \og champs de vecteurs unitaire\fg\ utilis\'ee dans cette d\'efinition. Supposons la relation d'\'equivalence ergodique et fixons-en une repr\'esentation hilbertienne $\pi$ sur un champ $H$ de base $X$. Si $\xi : X \to H$ est un champ de vecteurs invariant, l'application $x\mapsto \norm{\xi_x}_x$ est constante. Par suite, si $\xi^n$ est une suite presque invariante de champs de vecteurs, on peut esp\'erer obtenir un contr\^ole sur la variation en $x$ des applications $x\mapsto \norm{\xi^n_x}_x$, du moins  pour $n$ grand, ce qui permettrait en particulier affaiblir l'hypoth\`ese $\norm{\xi^n_x}_x = 1$ presque s\^urement tout en conservant une d\'efinition identique de propri\'et\'e T de Kazhdan. 

\bigskip

L'un des r\'esultats techniques importants de cet article montre que si l'on suppose qu'il n'y a pas de perte de masse dans la suite presque invariante $\xi^n$, alors la repr\'esentation $\pi$ contient une suite presque invariante unitaire au sens ci-dessus. Plus pr\'ecis\'ement nous d\'emontrons l'\'equivalence de la d\'efinition pr\'ec\'edente et de la d\'efinition suivante.

\bigskip

\begin{tabular}{m{.7cm}m{12.5cm}}
(T) & Toute repr\'esentation hilbertienne poss\'edant une suite $\xi^n$  de champs de vecteurs presque invariante domin\'ee non triviale, au sens o\`u $\norm{\xi^n}_1 = 1$ et $\norm{\xi^n_x}\leqslant g(x)$ pour une fonction $g\in L^1(X)$, poss\`ede un champ de vecteurs invariant non trivial.
\end{tabular}

\bigskip

La d\'emonstration de ce r\'esultat fait l'objet du paragraphe  \ref{pinv}  (voir aussi le paragraphe  \ref{Kazhdan} pour le cas non ergodique). 

\bigskip

Une \'elaboration des techniques utilis\'ees pour obtenir cette caract\'erisation conduit alors au th\'eor\`eme suivant, qui relie la notion d'ergodicit\'e forte au ph\'enom\`ene de concentration de la mesure, en termes de fonctions 1-lipschitziennes (cf. \textsection\ref{ptes}). (Nous avions d\'ej\`a constat\'e des liens  entre ergodicit\'e forte et concentration lors de la premi\`ere partie de cet article \cite[thm. 2]{Pichot04_I}.)

\bigskip

\begin{dfn} Soient $H$ un espace de Hilbert et $(\mu_n)_n$ une suite de mesures de probabilit\'es sur $H$. On suppose  les premiers moments
\[
m_1(\mu_n)=\int_{H} \norm{y}d\mu_n(y)\leqslant C<\infty
\]
unifom\'ement born\'es par une constante $C>0$. Nous dirons que  la suite d'espaces m\'etriques-mesur\'es $(H,\norm\cdot,\mu_n)$ forme une \textnormal{famille de Levy} si pour tout $\varepsilon>0$ on a,
\[
\inf_{f} \mu_n\{\abs{f-m}\leqslant \varepsilon\} \to_n 1,
\]  
o\`u l'infimum est pris sur les fonctions 1-lipschitziennes  $f : H\to \RI$  et  $m = \int_{H} fd\mu_n$ est la valeur moyenne de $f$.
\end{dfn}

\bigskip

\begin{thm}\label{Levy-intr} Soit $H$ un espace de Hilbert. Soit $\R$ une relation d'\'equivalence ergodique sur un espace de probabilit\'e $(X,\mu)$. Alors $\R$ est fortement ergodique si et seulement si pour toute suite presque invariante $\xi^n : X \to H$ pour la repr\'esentation triviale de $\R$ dans $H$, domin\'ee par une fonction $g\in L^1$, la famille $(H, \norm \cdot, \mu_n)$ est une famille de Levy, o\`u $\mu_n = \xi^n_*\mu$ est la pouss\'ee en avant de $\mu$ sur $H$ par le champ $\xi^n$.
\end{thm}

\bigskip

Nous d\'emontrons ensuite le th\'eor\`eme suivant, dont l'analogue pour les groupes est bien connu, cf. \cite{HarpeValette89}. Sa d\'emonstration repose notamment sur la caract\'erisation de la propri\'et\'e T d\'ecrite ci-dessus et sur le th\'eor\`eme de concentration pr\'ec\'edent (thm. \ref{Levy-intr}).

\bigskip

\begin{thm}\label{prox-intr}
Soit $\R$ une relation d'\'equivalence ergodique ayant la propri\'et\'e T. Soit $\pi$ une repr\'esentation  de $\R$ poss\'edant une suite $\xi^n$ presque invariante telle que $\norm{\xi^n_x}_x\leqslant g(x)$ pour une fonction $g\in L^1(X)$. Il existe quitte \`a extraire une suite $\zeta_n$ de champs invariants domin\'es par $g$ tels que 
\[
\norm {\xi_x^n -  \zeta^n_x}_x \to 0
\]
pour presque tout $x\in X$.
\end{thm}

\bigskip

Ce th\'eor\`eme permet par exemple de donner une preuve directe du fait que, si $\R=\R_\alpha$  est une relation d'\'equivalence obtenue par action libre ergodique pr\'eservant une mesure de probabilit\'e d'un groupe $\G$, alors $\G$ poss\`ede la propri\'et\'e T si et seulement si $\R$ la poss\`ede (cf. \textsection\ref{ptes} pour des r\'ef\'erences concernant ce r\'esultat).

\bigskip
\bigskip

6. Nous avons \`a pr\'esent \'evoqu\'e tous les ingr\'edients n\'ecessaires pour \'etudier les aspects spectraux de la propri\'et\'e T et ainsi d\'emontrer le th\'eor\`eme \ref{T-intr}. Nous en donnons la preuve au paragraphe \ref{dem}, qui contient \'egalement la preuve du th\'eor\`eme \ref{ef-intr}. Les d\'emonstrations utilisent l'ergodicit\'e forte comme outil essentiel et ceci a motiv\'e l'\'ecriture de la premi\`ere partie \cite{Pichot04_I} de cet article.  Nous renvoyons \`a cette premi\`ere partie, ainsi qu'\`a ses r\'ef\'erences, pour des rappels sur cette notion (notamment l'article de Hjorth-Kechris \cite{HjorthKechris03}).

\bigskip
\bigskip
\centerline{---}
\bigskip

Concluons par une application des id\'ees  pr\'ec\'edentes. Rappelons tout d'abord un crit\`ere bien connu, extrait de \cite{Zuk96}, pour qu'un groupe discret ait la propri\'et\'e de Kazhdan. On consid\`ere un groupe d\'enombrable de type fini $\G$, et on note $L$ le \og link\fg\ en l'identit\'e de son complexe de Cayley $Y$, de dimension 2, associ\'e \`a un syst\`eme g\'en\'erateur sym\'etrique fini ($L$ est le graphe fini donn\'e par l'intersection de $Y$ avec une sph\`ere de rayon suffisamment petit centr\'ee en l'identit\'e). On suppose $L$ connexe.

\bigskip

{\bf Crit\`ere $\lambda_1 > 1/2$.} \textit{Si $\lambda_1(L) >1/2$, alors $\G$ poss\`ede la propri\'et\'e T de Kazhdan.}

\bigskip

L'hypoth\`ese $\lambda_1 >1/2$ peut \^etre consid\'er\'ee comme une hypoth\`ese de \og courbure positive\fg\ (cf. \cite{Garland73}), et signifie par d\'efinition que la premi\`ere valeur propre non nulle du laplacien discret sur $L$ est strictement sup\'erieure \`a $1/2$.

\bigskip

Ce crit\`ere tire ses origines dans les travaux de Garland \cite{Garland73} sur l'annulation de la cohomologie (en certains degr\'es) de groupes d'automorphismes d'immeubles de Bruhat-Tits cocompacts. Gromov en a donn\'e une preuve nouvelle  dans  \cite{Gromov03}, bas\'ee sur l'\'etude des marches al\'eatoires sur les groupes discrets. Poursuivant ces id\'ees, nous d\'emontrons le r\'esultat suivant au paragraphe \ref{Garland}.

\bigskip

\begin{thm} Soit $\Q$ un espace singulier de type fini, $\tilde \Sigma$ un complexe simplicial $\Q$-p\'eriodique de dimension 2 uniform\'ement localement fini. On fixe une mesure quasi-invariante $\mu$ sur la lamination $X$ des sommets associ\'ee \`a $\tilde \Sigma$ et on consid\`ere un nombre r\'eel $\delta_\mu \geqslant 1$ tel que 
\[
\delta_\mu^{-1} \leqslant \delta(y,z) \leqslant \delta_\mu
\]
pour presque toute ar\^ete $(y,z)$ de $\tilde \Sigma$, o\`u $\delta$ est le cocyle de Radon-Nikodym associ\'e \`a $\mu$. On suppose que presque tout link $L$ de $\tilde \Sigma$ est connexe et v\'erifie $\lambda_1(L) \geqslant \lambda$, o\`u $\lambda$ est un nombre r\'eel v\'erifiant
\[
\lambda > \delta_\mu^3/2.
\]
Alors $\Q$ poss\`ede la propri\'et\'e T de Kazhdan.
\end{thm}

\bigskip

L'objet de ma note aux comptes-rendus \cite{Pichot03} \'etait d'observer que le th\'eor\`eme ci-dessus est vrai pour les relations d'\'equivalence de type $\IIi$, i.e. de traiter le cas $\delta_\mu =1$ avec les notations du th\'eor\`eme.

\bigskip
\bigskip

\centerline{---}

\bigskip

Je remercie Damien Gaboriau pour son aide constante au cours de l'\'elaboration de ce travail. 

Je dois \'egalement beaucoup \`a \'Etienne Ghys, ainsi qu'aux excellentes conditions de travail dont on b\'en\'eficie au sein de l'UMPA.

\bigskip
\bigskip

\section{Notations}

Les deux parties de cet article (le pr\'esent article, et \cite{Pichot04_I}) ont \'et\'e \'ecrites simultan\'ement. Pour cette raison il ne nous semble pas indispensable de r\'e\'ecrire ici une pr\'esentation g\'en\'erale du contexte qui nous concerne ; nous renvoyons \`a la premi\`ere partie \cite{Pichot04_I} et aux r\'ef\'erences pour cela. Contentons nous de rappeler les notations que nous utiliserons dans la suite.

\bigskip

\'Etant donn\'e un  espace bor\'elien standard  $X$ muni d'une mesure de probabilit\'e sans atome $\mu$, on dit que le couple $(X,\mu)$ est un espace de probabilit\'e.
Soit $(X,\mu)$ un espace de probabilit\'e.  Une relation d'\'equivalence  \`a classes d\'enombrables $\R$ sur $(X,\mu)$ est mesur\'ee si son graphe $\R\subset X\times X$ est bor\'elien et si $\mu$ est une mesure quasi-invariante, au sens o\`u le satur\'e d'un bor\'elien n\'egligeable est n\'egligeable. On note $\h$ la mesure sur $\R$ associ\'ee \`a $\mu$ et au syst\`eme de Haar canonique $(\h^x)_{x\in X}$ de d\'ecompte horizontal.  Explicitement,
\[
\h(K) = \int_X \h^x(K) d\mu(x)=\int_X \#K^x d\mu(x)
\]
o\`u $K\subset \R$ est une partie bor\'elienne et $K^x = \{(x,y)\in K\cap \R\}$.  Le symbole $\h_1$ d\'esigne une mesure de probabilit\'e sur $\R$ \'equivalente \`a $\h$. Les espaces singuliers sont not\'es $\Q$ (e.g. $\Q=X/\R$) et suppos\'es munis d'une mesure transverse invariante $\Lambda$.  Les lettres $A$, $B$,... sont des parties bor\'eliennes de $X$. La lettre $K$ d\'esigne une partie mesurable sym\'etrique de $\R$ et est consid\'er\'ee comme une structure simpliciale mesurable sur les classes (d'une sous-relation) de $\R$. La lettre $\Sigma$ est attach\'ee aux  complexes simpliciaux $\Q$-p\'eriodiques. La lettre $\G$ d\'esigne un groupe d\'enombrable. 

\bigskip

On notera $H$ les champs mesurables d'espaces de Hilbert s\'eparables de base $X$, et $\pi$ les repr\'esentations unitaires de $\R$ sur $H$. Les champs de vecteurs sur $X$, i.e. les sections de $H$, sont not\'es $\xi$, $\zeta$,...  Rappelons qu'une repr\'esentation unitaire de $\R$ sur $H$ consiste en la donn\'ee d'une famille d'op\'erateurs unitaires
\[
\pi(x,y) : H_y \to H_x,
\]
$(x,y)\in \R$, satisfaisant aux conditions de composition et de mesurabilit\'e suivantes :
\begin{itemize}
\item $\pi(x,x)=\id\ \mathrm{et} \  \pi(x,z)= \pi(x,y)\pi(y,z)$ pour tout $x\sim y \sim z$.
\item les coefficients
\[
(x,y) \mapsto \langle\pi(x,y) \xi_y |  \eta_x \rangle_x
\]
sont mesurables pour tous champs de vecteurs mesurables $\xi, \eta : X \to H$.
\end{itemize}

\bigskip
 Par exemple, la repr\'esentation r\'eguli\`ere de $\R$ sur  le champ d'espaces de Hilbert 
\[
H : x\mapsto \ell^2(\R^x,\h^x)
\]
qui \`a tout point $x\in X$ associe l'espace des fonctions de carr\'e int\'egrable sur la  classe d'\'equivalence de $x$ (pour la mesure de d\'ecompte horizontal $\h^x$), 
\[
\pi(x,y) : \ell^2(\R^y) \to \ell^2(\R^x)
\] 
est d\'efinie par $\pi(x,y)f(x,z) = f(y,z)$.\\

Nous utiliserons la lettre  $\nu$ pour les marches al\'eatoires sur  un graphage de relation d'equivalence ($\tilde \nu$ pour les marches al\'eatoires sur le graphe quasi-p\'eriodique associ\'e). 

\bigskip

Enfin, la lettre $D$ d\'esigne une contraction hilbertienne (e.g. diffusion hilbertienne associ\'ee \`a $\nu$), et $E$ est la fonction \'energie associ\'ee. On distingue les \'energies $E_n$ associ\'ees aux diffusions $D^n$. Les diverses \og constantes de diffusion\fg\ sont not\'ees $\kappa$, $\lambda$ et $c_n$.

\bigskip
\bigskip
\centerline{---}
\bigskip

Nous utiliserons le r\'esultat suivant, extrait de \cite{Pichot04_I}.

\begin{dfn} On dit que $K$ poss\`ede des \textnormal{suites de F\o lner \'evanescentes}  relativement \`a $\mu$ s'il existe une suite $(A_n)$ de bor\'eliens non n\'egligeables de $X$ et une suite $(\varepsilon_n)$ de nombres r\'eels convergeant vers 0  telles que  
\[
\mu(A_n) \to 0 \qquad \mathrm{et} \qquad \mu(\del_{K} A_n) \leqslant \varepsilon_n \mu(A_n).
\]
\end{dfn}

\bigskip

\begin{thm}
Soit $\R$ une relation d'\'equivalence ergodique de type fini pr\'eservant une mesure de probabilit\'e $\mu$. Alors $\R$ poss\`ede un quotient moyennable si et seulement si chacun de ses graphages u.l.f. contient des suites de F\o lner \'evanescentes.
\end{thm}

\bigskip
\bigskip

\vspace{1cm}
\section{Un lemme technique}\label{pinv}

Soit $\R$ une relation d'\'equivalence mesur\'ee sur un espace de probabilit\'e $(X,\mu)$ et $\pi$ une repr\'esentation de $\R$ sur un champ mesurable d'espaces hilbertiens $H$ de base $X$.

\medskip

Dans ce paragraphe nous \'etudions le probl\`eme de dispersion de masse des champs presque invariants de $\pi$.

\medskip

\begin{dfn} Soient une partie bor\'elienne $K \subset \R$ et un nombre r\'eel strictement positif $\varepsilon$.  On dit qu'un champ de vecteurs $\xi : X \to H$ (i.e. une section mesurable de $H$) est \textnormal{$(K,\varepsilon)$-invariant} si 
\[
\norm{\pi(x,y)\xi_y - \xi_x} \leqslant \varepsilon, 
\]
pour $\h$-presque tout  $(x,y) \in K$.
\end{dfn}

\bigskip

\begin{lem} Soit $\h_1$ une mesure de probabilit\'e sur $\R$ \'equivalente \`a $\h$. 
Soit $(K_n)_n$ une suite de parties bor\'eliennes de $\R$. Si $\h_1(K_n) \to 1$, il existe une suite extraite $K_{m_1},K_{m_2},\ldots $ telle que  la suite $(K'_i)_{i\geqslant 1}$ d\'efinie par $K'_i = \cap_{m_j\geqslant m_i}K_{m_j}$ soit une approximation croissante de $\R$ (au sens o\`u $\R = \cup_i K'_i$ \`a un n\'egligeable pr\`es). Si inversement $\R = \cup K_n$ est une approximation croissante de $\R$, alors $\h_1(K_n) \to 1$.
\end{lem}

\begin{dem}
Notons  $A_n=\R \backslash K_n$. Par hypoth\`ese $\h_1(A_n) \to 0$ et, quitte \`a extraire, on peut supposer $\sum \h_1(A_n) < \infty$. Le lemme de Borel-Cantelli montre que $\h_1(\limsup A_n) = 0$, o\`u $\limsup A_n = \cap_i\cup_{n\geqslant i} A_n$. Ainsi, en notant $K_i'=\cap_{n\geqslant i} K_n$, on obtient $\h_1(K_i') \to 1$, d'o\`u le r\'esultat.

R\'eciproquement consid\'erons la suite de fonctions indicatrices $\chi_{K_n} \in L^\infty(\R)$. Par hypoth\`ese elle converge presque s\^urement vers la fonction constante \'egale \`a 1 sur $\R$, et la convergence a lieu dans $L^1(\R,\h_1)$ \'egalement.
\end{dem}

\begin{dfn}
 On dit qu'une suite $\xi^n :  X \to H$ de champs de vecteurs $(K_n,\varepsilon_n)$-invariants est \textnormal{presque invariante} si $\epsilon_n$ converge vers 0 et si $K_n \subset \R$ est croissante et exhaustive. 
\end{dfn}

Le lemme pr\'ec\'edent montre qu'une suite de champs $(K_n,\varepsilon_n)$-invariants telle que $\varepsilon_n \to 0$ et $\h_1(K_n) \to 1$ pour une mesure de probabilit\'e $\h_1$ sur $\R$ \'equivalente \`a $\h$, est, quitte \`a extraire, une suite presque invariante au sens ci-dessus.\\

\textit{Exemple.} Il est facile de voir que toute repr\'esentation contient des suites presque invariantes $\xi^n : X \to H$ de carr\'e int\'egrable et de norme $\norm {\xi^n}_2=1$ (en choisissant par exemple $\xi^n=\chi_{A_n}\cdot\xi/\sqrt{\mu(A_n)}$ o\`u $\mu(A_n)\to 0$ et $\xi$ est un champ unitaire fixe).

\bigskip

 Le r\'esultat qui suit montre que toute repr\'esentation contenant des suites presque invariantes \textit{domin\'ees} de norme 1 contient des sections presque invariantes du fibr\'e en sph\`eres unit\'es de $H$.

\bigskip

\begin{lem}\label{lemtech}
Soit $\R$ une relation d'\'equivalence mesur\'ee ergodique sur un espace de probabilit\'e $(X,\mu)$. On fixe une repr\'esentation $\pi$ de $\R$ sur un champ hilbertien $H$ de base $X$. Les propri\'et\'es suivantes sont \'equivalentes.
\begin{itemize}
\item[i.] \cite{Zimmer84} Pour tout groupe d\'enombrable $\G$ et toute action $\alpha$ de $\G$ telle que  $\R = \R_\alpha$, pour tout $\varepsilon > 0$ et toute partie finie $K$ de $\G$, il existe $\xi\in L^\infty(X,H)$ tel que $\norm \xi_\infty = 1$ et 
\[
\mu\{ \abs{\langle\pi(x,\alpha(\gamma) x)\xi_{\gamma x}|\xi_x\rangle -1} \geqslant \varepsilon \} \leqslant \varepsilon
\]
pour tout $\gamma \in K$.

\item[ii.] \cite{Moore82} Il existe une suite presque invariante $\xi^n : X \to H$ de champs de vecteurs tels que  $\norm {\xi^n_x}_x = 1$ pour presque tout $x \in X$.

\item[iii.] [$L^\infty$] Il existe une suite presque invariante $\xi^n\in L^\infty(X,H)$ de champs de vecteurs et deux constantes $\eta \geqslant 1$ et $\delta >0$ telles que pour tout $n$
\[
\mu\{ \frac 1 \eta \leqslant \norm {\xi^n_x}_x \leqslant \eta\} \geqslant \delta.
\]

\item[iv.] [$L^p$, $1\leqslant p< \infty$] Il existe une suite presque invariante $\xi^n\in L^p(X,H)$ de champs de vecteurs tels que $\norm {\xi^n}_p = 1$ et une fonction positive $g \in L^p(X)$ telle que $\norm {\xi^n_x}_x \leqslant g(x)$.
\end{itemize}
\end{lem}

\bigskip
La fin du paragraphe est consacr\'ee \`a la d\'emonstration de ce r\'esultat.\\

Consid\'erons une relation d'\'equivalence ergodique $\R$ sur un espace de probabilit\'e $(X,\mu)$ et $\pi$ une repr\'esentation de $\R$ sur un champ hilbertien $H$ de base X.\\

\textit{Preuve de $i \impl ii$.} Soit $\G$ un groupe discret et $\R=\R_\alpha$ la relation d'\'equivalence associ\'ee \`a une action $\alpha$ de $\G$ sur $(X,\mu)$. Par hypoth\`ese, \'etant donn\'ee une exhaustion $K_n$  de $\G$ par parties finies, il existe une suite $\xi^n$ telle que $\norm {\xi^n}_\infty = 1$, et une suite $A_n$ de bor\'eliens telle que $\mu( A_n) \geqslant 1-1/n$, v\'erifiant, 
\[
\abs{\langle\pi(x,\gamma x)\xi^n_{\gamma x}|\xi^n_x\rangle -1} \leqslant 1/n
\]
presque s\^urement sur $A_n$ (pour tout $\gamma \in K_n$). En particulier $\norm{\xi^n_x}^2 \geqslant 1- 1/n$ sur $A_n$ ($e \in K_n$ pour $n$ suffisamment grand). Soit $F_n= \graph(K_n) \cap  A_n \times A_n$. On pose $\tilde \xi^n_x=\xi^n_x/\norm{\xi^n_x}$ sur $A_n$ et $\eta_x \in H_x$ quelconque de norme 1 ailleurs (mesurable). Alors $\h_1(F_n) \to 1$ et, pour $(x,y) \in F_n$, on a 
\begin{eqnarray}
\norm{\pi(x,y)\tilde \xi^n_y-\tilde \xi^n_x}^2 
&\leqslant& 2 \abs{1 -  {\langle\pi(x,y) \xi^n_y| \xi^n_x\rangle \over \norm {\xi^n_x}\norm {\xi^n_y}}}  \nonumber \\
&\leqslant& 2\cdot{1/n +1/n \over 1-1/n} \to 0. \nonumber
\end{eqnarray}

\bigskip

\textit{Preuve de $ii \impl i$}. Soit $\G$ un groupe discret et $\R=\R_\alpha$ la relation d'\'equivalence associ\'ee \`a une action $\alpha$ de $\G$ sur $(X,\mu)$. Soit $K$ une partie finie de $\G$ et $\varepsilon >0$ ;  comme $\h_1(\graph(K) \backslash F_n)\to 0$, l'ensemble
\[
A_n=\{x\in X,\ \exists \gamma \in K,\ (x, \gamma x) \notin F_n \} \to_\mu 0.
\]
En effet $(pr_h)_*({\h_1}_{|\graph(\gamma)})$ et $\mu$ sont \'equivalentes pour tout $\gamma \in \G$, o\`u $pr_h(x,y)=x$ est la projection horizontale. Par hypoth\`ese sur $X\backslash A_n$, $\norm {\pi(x, \gamma x)\xi^n_{\gamma x} - \xi^n_x} \leqslant \varepsilon_n$ pour tout $\gamma \in K$. Notons $z_n(x,y) = \langle \pi(x,y) \xi^n_y| \xi^n_x\rangle \in \CI$, de sorte que
\[
\mu \{  x \in X, \ 2- 2 \Re z_n(x,\gamma x) \geqslant \varepsilon_n \}\leqslant \mu(A_n) \to 0,
\]
pour tout $\gamma \in K$. Comme $ \abs{ z_n(x,y)} \leqslant 1$ presque s\^urement sur $\R$, on a donc 
\[
\mu \{  x \in X, \abs {1 -  z_n(x,\gamma x)}\geqslant \varepsilon/2 \} \to 0,
\]
d'o\`u le r\'esultat.

\bigskip

\textit{Preuve de $iii \impl ii$.} Fixons deux constantes $\eta \geqslant 1$ et $\delta >0$ et supposons qu'il existe une suite $\xi^n\in L^\infty(X,H)$ de champs de vecteurs mesurables presque invariants  telle que
\[
\mu\{ \frac 1 \eta \leqslant \norm {\xi^n_x}_x \leqslant \eta\} \geqslant \delta.
\]

Consid\'erons l'espace $\cal E$  des suites de couples $(\xi^n,\alpha_n)_n$ form\'es d'un champ de vecteurs et d'un nombre r\'eel positif telles que :

- il existe une suite $K_n$ de $\R$ et une suite d\'ecroissante $\varepsilon_n$ de nombre r\'eels positifs telles que $\xi^n$ soit $(K_n,\varepsilon_n)$-invariante, avec $\h_1(K_n) \to 1$ et $\varepsilon_n \to 0$.

- $\alpha_n$ soit une suite d\'ecroissante de nombres r\'eels tendant vers 0.

- $\mu(A_n)$  soit une suite convergente, o\`u
\[
A_n = \{  \frac 1 \eta - \alpha_n \leqslant \norm {\xi^n_x}_x \leqslant \eta + \alpha_n\}.
\]
On consid\`ere l'application $\cal E \to \RI$ d\'efinie par 
\[
\Psi : (\xi^n,\alpha_n)_n \mapsto \lim \mu(A_n),
\]
et on note $\overline \delta = \sup \Psi( {\mathcal E}) \leqslant 1$. Par hypoth\`ese $\overline \delta >0$. Montrons que $\overline \delta =1$.\\

Il existe par proc\'ed\'e diagonal un \'el\'ement de $\cal E$, disons $(\xi^n,\alpha_n)_n$, tel que $\mu(A_n) \to \overline \delta$. Notons $C_n = ((A_n \times X \cup X \times A_n) \backslash A_n \times A_n) \cap \R$.\\

Supposons qu'il existe une suite extraite $(\xi^{n_i},\alpha_{n_i})$ telle que $\h_1(C_{n_i}) \geqslant c$ pour un nombre r\'eel $c >0$. Alors il existe quitte \`a extraire un nombre $c'>0$ tel que $\mu(pr_h(C_{n_i}\cap X \times A_{n_i})) \geqslant c'$ ou $\mu(pr_v(C_{n_i}\cap A_{n_i} \times X)) \geqslant c'$. Rappelons que $pr_h : \R \to X$ est la projection horizontale d\'efinie par $pr_h(x,y)=x$ (et $pr_v(x,y)=y$ la projection verticale). Notons que (quitte \`a extraire une seconde fois) la suite $\sigma =(\xi^{n_i},\alpha_{n_i} + \varepsilon_{n_i})_{n_i}$ appartient \`a $\cal E$. Or 
\[
A_{n_i}^\sigma = \{\frac 1 \eta - \alpha_{n_i}-\varepsilon_{n_i} \leqslant \norm {\xi^{n_i}_x}_x \leqslant \eta + \alpha_{n_i} + \varepsilon_{n_i}\}
\]
contient $pr_h((C_{n_i}\cap X\times A_{n_i})\cap K_{n_i})$ et $pr_v((C_{n_i}\cap A_{n_i}\times X)\cap K_{n_i})$. Comme $\h_1(K_{n_i}) \to 1$, on a $\h_1(C_{n_i}\backslash K_{n_i}) \to 0$ donc  $\mu(pr_h((C_{n_i}\cap X\times A_{n_i})\cap K_{n_i})) \geqslant c'/2$ ou $\mu(pr_v((C_{n_i}\cap A_{n_i}\times X)\cap K_{n_i})) \geqslant c'/2$ pour ${n_i}$ grand. Ceci contredit la maximalit\'e de $\overline \delta$. \\

Par suite $\h_1(C_n) \to  0$.\\

Soit $\varphi \in [\R]$ un isomorphisme de $\R$. Notons que $(pr_h)_*({\h_1}_{|\graph(\varphi^{-1})}) \sim \mu$. Donc
\[
\mu(\varphi A_n\backslash A_n) = \mu\{x\in \varphi A_n \mid x \notin A_n \} \sim \h_1(\{(\varphi x,x),\, x\in A_n\}\cap C_n) \leqslant \h_1(C_n) \to 0.
\]
(On dit que $(A_n)$ est asymptotiquement invariante.) Il en r\'esulte que toute limite faible de la suite $\chi_{A_n} \in L^\infty(X)$ est constante (par ergodicit\'e), et on a donc  
\[
\lim \mu(A_n \cap C) - \mu(A_n) \mu(C) = 0
\] 
pour tout bor\'elien $C \subset X$ (\cite[page 95]{ConnesKrieger77,JonesSchmidt87}). Par suite pour tout $n$ on a
\[
\lim_m \mu(A_n \cap A_m) = \mu(A_n)\overline \delta >0.
\]
Supposons $\overline \delta < 1$ et consid\'erons, \'etant donn\'e $n$, un entier $m=m(n)>n$ suffisament grand pour  que 
\[
\abs{\mu(A_n \cap A_m) -\mu(A_n)\overline \delta} \leqslant \frac 1 n
\]
et
\[
\abs{\h_1(K_n \cap K_m) - \h_1(K_n)} \leqslant \frac 1 n.
\]
Construisons alors une suite $\sigma \in \cal E$ de la fa\c con suivante. Pour tout $n$ on pose $\xi^n_\sigma = \xi^n$ sur $A_n$ et $\xi^n_\sigma = \xi^{m(n)}$ sur $A_{m(n)}\backslash A_n$ (et 0 ailleurs). Alors $\sigma = (\xi^n_\sigma,\alpha_n) \in \cal E$. En effet il suffit de choisir $\varepsilon_n^\sigma = \varepsilon_n$ et $K_n^\sigma = K_n\cap K_{m(n)} \backslash C_n \cup C_{m(n)}$. Par ailleurs $\lim \mu(A_n^\sigma) = 2\overline \delta - \overline \delta^2 > \overline \delta$, d'o\`u une contradiction.\\

Finalement $\overline \delta = 1$. Consid\'erons alors la suite $\tilde \xi_n$  d\'efinie par $\tilde \xi^n_x=\xi^n_x/\norm{\xi^n_x}$ sur $A_n$ et $\tilde \xi_x \in H_x$ quelconque de norme 1 ailleurs (mesurable). Notons $\tilde K_n = K_n \cap (A_n\times A_n)$. Alors $\h_1(\tilde K_n) \to 1$, et $\tilde \xi^n$ est $(\tilde K_n, \eta\cdot \varepsilon_n)$-invariant.
\bigskip

\textit{Preuve de $iv \impl iii$.} Soit $p \in[1, \infty[$. Soit $\eta\geqslant 3$ suffisament grand pour que la norme de $g^p$ restreinte \`a  $A=\{ g > \eta\}$ soit $\leqslant 1/2$. Alors, en notant $A_n=\{ \frac 1 \eta \leqslant \norm {\xi^n_x}_x \leqslant \eta\}$, on a
\[
1 = \norm {\xi^n}_p^p = \int_{A_n} \norm {\xi^n}^p + \int_{X \backslash A_n} \norm {\xi^n}^p \leqslant \eta^p \mu(A_n) + 1/2 + \frac 1 {\eta^p} \mu(X\backslash A_n)
\]
donc
\[
\mu(A_n) \geqslant \frac{\eta^p-2}{2\eta^{2p}- 2 }.
\]
Consid\'erons alors une suite $c_n$ de nombres r\'eels $\geqslant \eta$ de sorte que $\mu(C_n) \to 1$ o\`u $C_n = \{\norm{\xi_n} \leqslant c_n\}$. Soit $\tilde \xi^n = \xi^n_{|C_n}$ et $\tilde K_n = K_n \cap C_n\times C_n$. Alors $\tilde \xi_n \in L^\infty$ est une suite $(\tilde K_n, \varepsilon_n)$-invariante satisfaisant aux hypoth\`eses de \textit{iii}.

\bigskip

Enfin $ii \impl iv$ est trivial.

\bigskip

\textit{Remarque.} Notons que si $\R$ n'est pas ergodique, \textit{ii} et \textit{iii} ne sont pas \'equivalentes,  comme le montre l'exemple \'el\'ementaire d'une relation d'\'equivalence obtenue comme r\'eunion disjointe de deux relations ergodiques, dont l'une est  munie d'une repr\'esentation sans champ presque invariant et l'autre d'une repr\'esentation ayant des champs presque invariants unitaires au sens du lemme.

\bigskip

Le lemme suivant r\'eunit quelques faits g\'en\'eraux utilis\'es au cours de la d\'emonstration ci-dessus.

\bigskip

\begin{lem} Soit $\R$ une relation d'\'equivalence mesur\'ee sur un espace de probabilit\'e $(X,\mu)$ et $\h$ la mesure de d\'ecompte horizontal sur $\R$ associ\'ee \`a $\mu$. Soit $\G$ un groupe d\'enombrable et $\alpha$ une action de $\G$ telle que $\R=\R_\alpha$.
\begin{itemize}
\item Si $\h_1, \h_1'$ sont deux mesures de probabilit\'e sur $\R$ \'equivalentes \`a $\h$, alors
\[
\h_1(K_n) \to 1 \ssi \h_1'(K_n) \to 1,
\]
o\`u $(K_n)$ est une suite de parties bor\'eliennes de $\R$.
\item Soit $(F_n)$ une suite croissante de parties de $\G$. Soit $K_n = \alpha(F_n) \subset \R$ la r\'eunion des graphes des \'el\'ements de $F_n$. Alors $(F_n)$ est une suite exhaustive de  $\G$ si et seulement si $\h_1(K_n) \to 1$.
\item Soient $(A_n),\, (B_n)$ deux suites de parties bor\'eliennes de $X$. Alors
\[
\mu(A_n) \to 1\ \mathrm{et}\ \mu(B_n)\to 1 \ssi \h_1(A_n \times B_n\cap \R) \to 1.
\]
\item Soit $\h_1$ une mesure de probabilit\'e sur $\R$ \'equivalente \`a $\h$. Soit $F\subset \R$ une partie bor\'elienne telle que $\h(F) < \infty$.  On consid\`ere une suite $F_n\subset \R$ de parties bor\'eliennes telle que $\h_1(F_n) \to 1$. Alors $\h(F \backslash F_n) \to 0$.
\end{itemize}
\end{lem}

\begin{dem} Rappelons que si $\h$ et $\h'$ sont deux mesures finies telles que $\h \ll \h'$ sur un espace bor\'elien standard, alors $\h'(A_n) \to 0 \impl \h(A_n) \to 0$ pour toute suite de bor\'eliens $A_n$. Ceci d\'emontre la premi\`ere affirmation. La seconde est triviale. Pour la troisi\`eme, notons que $(p_v)_*(\h_1)$ et $\mu$ sont \'equivalentes, donc $\mu(A_n) \to 1 \ssi \h_1(A_n \times X\cap \R) \to 1$. Or $\h_1(A_n \times  X \cap X\times B_n\cap \R) \to 1 \ssi  \h_1(A_n \times  X\cap \R) \to 1$ et  $\h_1(X \times B_n\cap \R) \to 1$. Pour la derni\`ere affirmation on a $\h_1(\R \backslash F_n) \to 0$, donc  $\h_1(F \backslash F_n) \to 0$. Comme $\h_{|F} \ll \h_1$, on a \'egalement $\h(F \backslash F_n) \to 0$.
\end{dem}

\bigskip

Terminons enfin par une variation, pour  r\'ef\'erences ult\'erieures, sur le  d\'ebut de l'implication $iii\impl ii$ du lemme \ref{lemtech}.

\bigskip

\begin{lem}\label{pi-ai}
Soit $\R$ une relation d'\'equivalence mesur\'ee sur un espace de probabilit\'e $(X,\mu)$. On fixe une repr\'esentation $\pi$ de $\R$ sur un champ hilbertien $H$ de base $X$. \'Etant donn\'es un champ mesurable $C=(C_x)_{x\in X}$ de parties mesurables de $H$ et un nombre r\'eel positif $\alpha$, on note $C_\alpha$ le champ mesurable qui \`a $x\in X$ associe le $\alpha$-voisinage $(C_x)_\alpha$ de $C_x$ dans $H_x$. 

Soient $(\xi^n)_{n\geqslant 0}$ une suite presque invariante de champs de vecteurs et $(C^n)_{n\geqslant 0}$ une suite de champs mesurables et invariants (i.e. stable par $\pi$) de parties mesurables de  $H$.  Il existe une suite extraite $(\xi_{m_1}, \xi_{m_2}, \ldots)$ de $(\xi_n)$ et une suite d\'ecroissante $(\alpha_{m_1},\alpha_{m_2},\ldots)$ de nombres r\'eels convergeant vers 0, de sorte que la suite $(A_{m_i})$ d\'efinie par
\[
A_{m_i} = \{ \xi^{m_i}_x \in (C_x^{m_i})_{\alpha_{m_i}}\}
\]
soit asymptotiquement invariante.
\end{lem}

\begin{dem} Reprenons la demonstration du lemme \ref{lemtech}. Consid\'erons l'espace $\cal E$  des suites de couples $(\xi^m,\alpha_m)_m$ form\'es d'un champ de vecteurs et d'un nombre r\'eel positif telles que :

- $\xi^m$ est une suite extraite de la suite presque invariante $\xi^n$,

- $\alpha_m$ est une suite d\'ecroissante de nombres r\'eels tendant vers 0,

- $\mu(A_m)$  soit une suite convergente, o\`u
\[
A_{m} = \{ \xi^{m}_x \in (C_x^{m})_{\alpha_{m}}\}.
\]
On consid\`ere encore l'application $\Psi : \cal E \to \RI$ d\'efinie par $(\xi^m,\alpha_m)_m \mapsto \lim \mu(A_m)$, et on note $\overline \delta = \sup \Psi( {\mathcal E}) \leqslant 1$. Si $\overline \delta =0$ ou $\overline \delta =1$ le r\'esultat est clair. Sinon il existe par proc\'ed\'e diagonal un \'el\'ement de $\cal E$, disons $(\xi^m,\alpha_m)_m$, tel que $\mu(A_m) \to \overline \delta\in ]0,1[$. Notons $E_n = ((A_n \times X \cup X \times A_n) \backslash A_n \times A_n) \cap \R$. On montre de m\^eme que dans le lemme $\ref{lemtech}$ en consid\'erant la suite $\sigma =(\xi^{m},\alpha_{m} + \varepsilon_{m})_{m}$ que s'il existe une suite extraite $(\xi^{m_i},\alpha_{m_i})$ de $(\xi^m,\alpha_m)_m$ telle que $\h_1(E_{m_i}) \geqslant c$ pour un nombre r\'eel $c >0$ alors
$\limsup \mu(A_{m_i}^\sigma)> \overline \delta$ o\`u  
\[
A_{m_i}^\sigma = \{\xi^{m_i}_x \in (C_x^{m_i})_{\alpha_{m_i}+\varepsilon_{m_i}}\}.
\]
Par suite $\h_1(E_m) \to  0$ et il en r\'esulte que $A_m$ est asymptotiquement invariante. 
\end{dem}

\vspace{1cm}
\section{La propri\'et\'e T de Kazhdan}\label{Kazhdan}

La propri\'et\'e T pour les groupes localement compacts a \'et\'e introduite par Kazhdan  \cite{Kazhdan67} en 1967. Nous ne ferons ici que rappeler la d\'efinition donn\'ee par Kazhdan, renvoyant \`a \cite{HarpeValette89,Valette02,BekkaHarpeValette04} pour des d\'etails. Nous adaptons ensuite cette d\'efinition aux espaces mesur\'es singuliers.
 
\bigskip

\subsection{Le cas des groupes.} Soit $\G$ un groupe localement compact. \\

On dit qu'une repr\'esentation unitaire $\pi$ de $\G$ sur un espace de Hilbert $H$ \textit{poss\`ede presque des vecteurs invariants} s'il existe une suite $\xi_n$ de vecteurs de norme 1, une suite exhaustive croissante $S_n$ de parties compactes de $\G$, et une suite $\varepsilon_n$ de nombres r\'eels tendant vers 0, telles que
\[
\norm{\pi(s)\xi_n - \xi_n }\leqslant \varepsilon_n, \quad \forall s\in S_n.
\]
On dit qu'un vecteur $\xi$ est invariant si $\pi(s)\xi=\xi$ pour tout $s\in \G$.\\

\begin{dfn}[\cite{Kazhdan67}]
On dit que $\G$ \textnormal{a la propri\'et\'e T de Kazhdan} si toute repr\'esentation fortement continue ayant presque des vecteurs invariants a des vecteurs invariants non nuls.
\end{dfn}

\medskip

\textit{Exemples de groupes de Kazhdan.} Kazhdan \cite{Kazhdan67} a  montr\'e que $\SL_3(\ZI)$ poss\`ede la propri\'et\'e T, de m\^eme que les r\'eseaux d'un groupe de  Lie simple de rang r\'eel sup\'erieur \`a 2. Les r\'eseaux de $\SL_2(\RI)$ ne l'ont pas alors que les r\'eseaux de $\Sp(1,n)$ l'ont \cite{HarpeValette89,Gromov03}. 

 De plus, la propri\'et\'e T est g\'en\'erique dans \og l'adh\'erence\fg\ de groupes hyperboliques, cf. \cite{Champetier00}, ainsi que pour certains mod\`eles statistiques de groupes al\'eatoires, cf. \cite{Zuk03,Gromov03}. 

Enfin rappelons qu'il existe un crit\`ere g\'eom\'etrique local, portant sur la \og courbure p-adique\fg\  \cite{Garland73} d'un poly\`edre fini, qui permet d'en d\'eduire la propri\'et\'e T pour son groupe fondamental. Voici quelques r\'ef\'erences \`a ce propos : \cite{Garland73,Borel74,Zuk96,BallSwiat97,CMS94,Gromov03,Ghys03,Zuk03} --- dont l'article original de Garland et l'\'enonc\'e bien connu de \.Zuk concernant sp\'ecifiquement la propri\'et\'e T (que nous avons rappel\'e en introduction). Notons que ce crit\`ere s'\'etend aux actions propres cocompactes \cite{Skandalis03}, et s'applique ainsi \`a $\SL_3(\QI_p)$ dont le classifiant propre, son immeuble de Bruhat-Tits, est l'exemple canonique de poly\`edre satisfaisant \`a ce crit\`ere. (cf. \'egalement \textsection\ref{Garland})\\

\subsection{Le cas des relations d'\'equivalence.} L'existence de repr\'esentations hilbertiennes int\'eressantes pour les relations d'\'equivalence mesur\'ees s'accompagne d'une notion de propri\'et\'e T, qui fut introduite par Moore et Zimmer au d\'ebut des ann\'ees 1980 (cf.  \cite{Zimmer84,Moore82} et \textsection\ref{pinv}). La d\'efinition est la suivante.

\medskip

Soit $\R$ une relation d'\'equivalence mesur\'ee sur un espace de probabilit\'e $(X,\mu)$. On dit qu'une repr\'esentation $\pi$ de $\R$ sur un champ hilbertien $H$ \textit{contient presque des champs invariants} s'il existe une suite $\xi_n$ de champs de vecteurs tels que $\norm{\xi_x^n} = 1$ pour presque tout $x\in X$, une suite exhaustive croissante $(K_n)$ de $\R$, et une suite $(\varepsilon_n)$ de nombres r\'eels tendant vers 0, telles que
\[
\norm{\pi(x,y)\xi_y^n - \xi_x^n }\leqslant \varepsilon_n, \quad \forall (x,y)\in K_n.
\]
On dit qu'un champ de vecteurs $\xi$ est invariant si $\pi(x,y)\xi_y=\xi_x$ pour presque tout $(x,y)\in \R$.\\

\begin{dfn}
 On dit que $\R$ \textnormal{poss\`ede la propri\'et\'e T de Kazhdan} si toute repr\'esentation unitaire ayant presque des champs invariants poss\`ede un champ invariant $\xi$ tel que $\norm{\xi_x}=1$ pour presque tout $x\in X$.
\end{dfn}

\bigskip
Remarquons que cette propri\'et\'e ne d\'epend que de la classe de $\mu$.

\bigskip

\begin{dfn}
On dit qu'un espace singulier $\Q$ poss\`ede la propri\'et\'e T si toutes ses d\'esingularisations discr\`etes la poss\`edent.
\end{dfn}

\bigskip

\begin{prop} 
La propri\'et\'e T de Kazhdan est un invariant d'isomorphisme stable de relations d'\'equivalence mesur\'ees.
\end{prop}

\begin{dem} Soit $\R$ une relation d'\'equivalence mesur\'ee sur $X$ et $\S$ la relation obtenue par restriction \`a un bor\'elien $Y$ de $X$ rencontrant presque toutes les orbites de $\R$.

Supposons que $\S$ poss\`ede la propri\'et\'e T. Soit $\pi$ une repr\'esentation de $\R$ sur $H$ et $\xi^n$ une suite de champs unitaires presque invariants. On note $H_Y$ la restriction de $H$ \`a $Y$ et $\pi_\S$ la restriction de $\pi$ \`a $\S$ (agissant sur $H_Y$). Alors $(\xi_n)_{|Y}$ est une suite de champs unitaires $\pi_\S$-presque invariant. Comme $\S$ a la propri\'et\'e T, il existe un champ $\xi : Y \to H_Y$ unitaire invariant par $\pi_\S$.  D\'efinissons $\overline \xi : X \to H$ par 
\[
\overline \xi_x = \pi(x,y)\xi_y
\]
pour $(x,y)\in \R$ tel que $x \notin Y$ et $y\in Y$. Comme $\xi$ est $\pi_\S$-invariant, $\overline \xi$ est d\'efini sans ambiguit\'e. Il est $\pi$-invariant par d\'efinition.

R\'eciproquement soit $\pi$ une repr\'esentation de $\S$ sur un champ $H$ de base $Y$. Fixons une partition bor\'elienne $Y_1=Y, Y_2, \ldots$ de $X$  et une famille $\varphi_2,\varphi_3,\ldots$ d'isomorphismes partiels de $\R$ surjectifs $\varphi_i : Y \surj Y_i$. Notons $\varphi_1 : Y \to Y$ l'identit\'e de $Y$. On d\'efinit un champ hilbertien $\overline H$ de base $X$  en associant \`a tout $y \in Y_i$ l'espace de Hilbert $H_{\varphi_i^{-1}y}$ et l'expression
\[
\overline \pi(x,y) = \pi(\varphi_i^{-1}x,\varphi_j^{-1}y),
\]
o\`u $x \in Y_i$ et $y\in Y_j$, d\'efinit une repr\'esentation de $\R$ sur $\overline H$. Une section $\xi$ de $H$ s'\'etend \`a $\overline H$ en posant $\overline \xi_y = \xi_{\varphi_i^{-1}y}$, $y\in Y_i$. Soit $\xi^n$ une suite $(K_n,\varepsilon_n)$-invariante de champs unitaires sur $Y$, o\`u $(K_n)$ est une suite croissante exhaustive de $\S$ et $\varepsilon_n \to 0$. On obtient en notant  
\[
\overline K_n = \amalg_{i,j} \varphi_{i,j}(K_n),
\]
o\`u $\varphi_{i,j}(x,y)= (\varphi_i x,\varphi_j y)$, une suite croissante exhaustive de $\R$. De plus $\overline \xi^n$ est une suite $(\overline K_n,\varepsilon_n)$-invariante pour $\overline \pi$. Si alors $\R$ poss\`ede la propri\'et\'e T, il existe un champ invariant unitaire $\xi$, et sa restriction \`a $Y$ est un champ unitaire $\pi$-invariant.
\end{dem}

\bigskip

Ainsi la propri\'et\'e T est un exemple de \og param\`etre de quasi-p\'eriodicit\'e\fg, suivant le point de vue de  \cite{Pichot04_I}.

\bigskip
\subsection{Une caract\'erisation des espaces de Kazhdan.}\label{nonerg}

Les r\'esultats du paragraphe \ref{pinv} (Lemme \ref{lemtech}) montrent qu'il est possible pour les relations ergodiques de caract\'eriser la propri\'et\'e T en termes de champs presque invariants domin\'es. Dans ce paragraphe nous \'enon\c cons explicitement ce r\'esultat, en l'\'etendant aux relations d'\'equivalence non ergodiques.  (Voir \'egalement la remarque suivant le lemme \ref{lemtech}.)

\bigskip

 Soit $H$ un champ mesurable d'espace de Hilbert sur un espace de probabilit\'e $(X,\mu)$. On dit qu'un champ de vecteurs sur $X$ est \textit{\`a support total} si son support $\{\xi_x \neq 0\}$ est de mesure 1, et qu'il est \textit{non trivial} si son support est de mesure non nulle.

\bigskip

\begin{thm}\label{dom} Une relation d'\'equivalence mesur\'ee poss\`ede la propri\'et\'e T de Kazhdan si et seulement si chacune de ses repr\'esentations hilbertiennes qui poss\`ede une suite presque invariante de champs de vecteurs, non triviale au sens du lemme \ref{lemtech}, iii ou iv, contient des champs invariants non triviaux.
\end{thm}

\begin{dem}
Supposons d'abord que la conclusion soit vraie et montrons que $\R$ poss\`ede la propri\'et\'e T.  Soit $\pi$ une repr\'esentation contenant une suite $\xi^n$ presque invariante, non triviale au sens \textit{ii}. Elle l'est donc au sens \textit{iii} (ou \textit{iv}) et il existe un champ invariant non trivial $\xi$. Son support, disons $\Omega$, est invariant (et non n\'egligeable). En consid\'erant la repr\'esentation coincidant avec $\pi$ sur $X\backslash \Omega$ et indentiquement nulle sur $\Omega$, on obtient une nouvelle repr\'esentation de $\R$ contenant presque des champs invariants au sens $iii$ (ou $\textit{iv}$). Il existe donc un champ invariant qui prolonge $\xi$ \`a un bor\'elien non n\'egligeable de $X\backslash \Omega$ ; \`a l'aide du lemme de Zorn, on construit ainsi facilement un champ invariant pour $\pi$ \`a support total. \\

R\'eciproquement supposons que $\R$ ait la propri\'et\'e T. Reprenons la d\'emonstration $iii \impl ii$ du lemme \ref{lemtech}. Avec des notations identiques, on obtient pour  $\eta \geqslant 1$ une suite asymptotiquement invariante
\[
A_n = \{ \frac 1 \eta-\alpha_n \leqslant \norm {\xi^n} \leqslant \eta+\alpha_n \}
\]
associ\'ee \`a une suite $\xi^n$  presque invariante non triviale au sens \textit{iii}, telle que 
\[
\mu(A_n)\to_n \overline \delta>0.
\]
Soit 
\[
(\R,\mu,\h) = \int_Z (\R_z, \mu_z,\h^z) d\mu_Z(z)
\]
la d\'esint\'egration de $\R$ en composante ergodique \cite{FeldmanMoore77}. Ainsi pour $\mu_Z$-presque tout $z\in Z$, $\R_z$ est une relation d'\'equivalence $\mu_z$-ergodique sur un espace bor\'elien standard $Y$, et $\R$ est isomorphe \`a la relation $(z,y)\sim (z',y')$ si et seulement si $z=z'$ et $y\R_z y'$ sur $Z\times Y$. Par ergodicit\'e, toute limite faible de la suite $\chi_{A_n} \in L^\infty(X)$ est une fonction ne d\'ependant que de $z\in Z$. Consid\'erons une telle fonction, disons  $\tilde \delta : X \to [0,1]$, et supposons quitte \`a extraire que $\chi_{A_n}$ converge faiblement vers $\tilde \delta$. Notons $\Omega = \{z\in Z,\, \tilde \delta_z \neq 0\}$ (non n\'egligeable) et $\tilde A_n=A_n\cap \Omega\times Y$ (ainsi $\mu(A_n \backslash \tilde A_n) \to 0$). On a, pour tout $n$ fix\'e,
\[
\lim_m \mu(\tilde A_n \cap\tilde  A_m) =\int_{\tilde A_n}\tilde \delta_z d\mu_Z(z).
\]
Supposons qu'il existe un bor\'elien $\Omega'\subset \Omega$ non trivial sur lequel $\tilde \delta_z < 1$ et consid\'erons, \'etant donn\'e $n$, un entier $m=m(n)>n$ suffisament grand pour  que 
\[
\abs{\mu(\tilde A_n \cap \tilde A_m) -\int_{\tilde A_n}\tilde \delta_z d\mu_Z(z)} \leqslant \frac 1 n
\]
et
\[
\abs{\h_1(K_n \cap K_m) - \h_1(K_n)} \leqslant \frac 1 n.
\]
Construisons alors une suite $\sigma \in \cal E$ de la fa\c con suivante. \'Etant donn\'e $n$ on pose, pour tout $z\notin \Omega$, $\xi^n_\sigma(z,y) = 0$, pour tout $z\in \Omega\backslash \Omega'$, $\xi^n_\sigma(z) = \xi^n_z$, et pour tout $z\in \Omega'$, $\xi^n_\sigma(z) = \xi^n_z$ sur $\tilde A_n^z$ et $\xi^n_\sigma(z) = \xi^{m(n)}_z$ sur $\tilde A_{m(n)}^z\backslash \tilde A_n^z$ (et 0 sur $Y\backslash\tilde  A_n^z\cup\tilde  A_{m(n)}^z$). Alors $\sigma = (\xi^n_\sigma,\alpha_n) \in \cal E$. En effet il suffit de choisir $\varepsilon_n^\sigma = \varepsilon_n$ et $K_n^\sigma = K_n\cap K_{m(n)} \backslash C_n \cup C_{m(n)}$. Par ailleurs 
\[
\lim \mu(A_n^\sigma) = 2\overline \delta - \int_\Omega \tilde \delta^2_z d\mu_Z(z)> \overline \delta
\]
d'o\`u une contradiction. Par suite  $\tilde \delta = 1$  presque surement sur $\Omega$.\\

Consid\'erons alors la repr\'esentation $\tilde \pi$ de $\R$ coincidant avec $\pi$ sur $\Omega\times Y$ et \'egale \`a la repr\'esentation trivial sur le champ constant $\CI$ de base $Z\backslash \Omega \times Y \subset X$. Soit $\tilde \xi_n$ la suite d\'efinie par $\tilde \xi^i_x = 1$ si $x\in Z\backslash \Omega \times Y$,  $\tilde \xi^i_x=\xi^i_x/\norm{\xi^i_x}$ sur $\tilde A_n \subset \Omega\times Y$ et $\eta_x \in H_x$ quelconque de norme 1 ailleurs (mesurable). Notons $\tilde K_n = (K_n \cap (\tilde A_n\times \tilde A_n))\amalg \R_{|X\backslash \Omega}$. Alors $\h_1(\tilde K_n) \to 1$, et $\tilde \xi^i$ est un champ $(\tilde K_n, \eta\cdot \varepsilon_n)$-invariant au sens \textit{ii} pour $\tilde \pi$. Par suite $\tilde \pi$ admet un champ invariant \`a support total, et ce champ restreint \`a $\Omega \times Y$ est invariant non trivial pour $\pi$.
\end{dem}

\vspace{1cm}
\section{Quelques propri\'et\'es des espaces de Kazhdan}\label{ptes}

\subsection{Existence d'une mesure invariante.} Les espaces singuliers ayant la propri\'et\'e T de Kazhdan sont, comme l'a montr\'e Zimmer,  de type $\II$. Plus g\'en\'eralement, le premier groupe de cohomologie $H^1(\R,\RI)$ \`a coefficients r\'eels d'une relation d'\'equivalence ayant la propri\'et\'e T est trivial  (cf. \cite{Moore82,Zimmer84}). (On renvoie ici \`a \cite{Jolissaint00} et \cite{Anan04}.)

\bigskip

\begin{prop}[Zimmer]
Soit $\R$ une relation d'\'equivalence sur un espace de probabilit\'e $(X,\mu)$ ayant la propri\'et\'e T de Kazhdan. Il existe une mesure $\sigma$-finie sur $X$ \'equivalente \`a $\mu$ et invariante par $\R$.
\end{prop}

\bigskip

\begin{cor}
Un espace singulier ayant la propri\'et\'e T est de type $\II$.
\end{cor}

\bigskip

\subsection{Approximations.}

L'un des r\'esultats fondamentaux de l'article de Kazhdan \cite{Kazhdan67} est le fait que tout groupe d\'enombrable ayant la propri\'et\'e T est de type fini. (Ainsi les r\'eseaux d'un groupe de Lie de rang sup\'erieur sont de type fini.) Cet \'enonc\'e  a \'et\'e adapt\'e par Moore \cite{Moore82} aux relations d'\'equivalence mesur\'ees (modulo de l\'eg\`eres imperfections que nous rectifions ci-dessous).

\bigskip

\begin{prop}
Toute approximation croissante d'une relation d'\'equivalence $\R$ ayant la propri\'et\'e T est, \`a isomorphisme stable pr\`es, constante \`a partir d'un certain rang.
\end{prop}

\bigskip

Plus pr\'ecis\'ement, si $\R_n \subset \R$ est une suite croissante et exhaustive de sous-relations de $\R$, alors \textit{il existe un entier $N$ et un bor\'elien $\R_N$-invariant non trivial sur lequel $\R$ et $\R_N$ co\"\i ncident.}\\

\begin{dem} Soient $\R$ une relation d'\'equivalence ayant la propri\'et\'e T et $\R_n$ une approximation de $\R$ au sens ci-dessus. Consid\'erons pour tout $x \in X$ l'espace de Hilbert $H^n_x$ des fonctions de carr\'e sommable d\'efinies sur les $\R_n$-classes de l'orbite $\R.x$ (i.e. les fonctions $\xi : \R^x \to \CI$ telles que $\xi(x,y)=\xi(x,z)$ si $(y,z)\in R_n$ et $\sum_{\overline {(y,z)} \in \R^x/\R_n^x} \abs{\xi(y,z)}^2 < \infty$) et notons $H^n$ le champ mesurable associ\'e aux espaces $H^n_x$. L'action r\'eguli\`ere de $\R$ sur elle-m\^eme qui permute les fibres horizontales  $(y,x)(x,z)=(y,z)$ induit une repr\'esentation $\pi_n$ de $\R$ sur $H^n$. Le champ $\chi_n : x \mapsto \1_{\R_n^x} \in H_n$ est  invariant par $\pi_n(\R_n)$.

 Comme $\h_1(\R_n) \to 1$,  la repr\'esentation $\pi = \oplus \pi_n$ poss\`ede une suite presque invariante de champs de vecteurs unitaires. Soit $\xi=(\xi^1, \xi^2, \ldots)$ un champ invariant non trivial. L'une des composantes $\xi^N$ est non nulle et comme $\pi$ est diagonale, $\xi^N$ est un champ non trivial invariant par $\pi_N(\R)$, i.e. $\xi^N$ d\'efinit mesurablement une fonction $\xi^N_{\R.x}$ par $\R$-orbite. Par construction cette fonction est constante sur les $\R_n$-classes et, \'etant de carr\'e int\'egrable, elle atteint son maximum sur un nombre fini de ces classes. De plus, $(R_n)$ \'etant croissante, il existe un entier $k\geqslant N$ et un bor\'elien $\R$-invariant $\Omega \subset X$ non trivial sur lequel $\xi^N_{\R.x}$ atteint son maximum sur exactement une $\R_k$-classe de chaque $\R$-classe de $\Omega$.  Alors $\R$ et $\R_k$ co\"\i ncident sur le bor\'elien (non n\'egligeable) constitu\'e par ces $\R_k$-classes. 
\end{dem}

\bigskip

\begin{cor}
Toute approximation croissante de $\R$ par des relations ergodiques est constante \`a partir d'un certain rang.
\end{cor}

\bigskip

Ce corollaire avait  \'egalement \'et\'e obtenu par Sorin Popa \cite{Popa86}.

\bigskip

De m\^eme qu'un groupe de Kazhdan est de type fini, on a le r\'esultat suivant.

\begin{cor}
Un espace singulier ergodique ayant la propri\'et\'e T de Kazhdan est de type fini.
\end{cor}

\begin{dem} 
Soit $\Q$ un espace singulier ergodique ayant la propri\'et\'e T (qui est donc de type $\II$). Soit $p :X\to \Q$ une d\'esingularisation discr\`ete de type $\IIi$ de $\Q$ et $\R=\R_p$ la relation associ\'ee. Il est bien connu qu'il existe un isomophisme ergodique $\varphi\in [\R]$. Num\'erotons une partition de $\R$ en isomorphisme partiel et consid\'erons l'ensemble bor\'elien $K_n\subset \R$ constitu\'ee de $\varphi$ et des $n$ premiers isomorphismes partiels de cette partition. Soit $\R_n$ la relation engendr\'ee par $K_n$. \'Evidemment $(\R_n)$ exhauste $\R$. Donc $\R_n = \R$ \`a n\'egligeable pr\`es.
\end{dem}

\bigskip

Observons que tous les exemples connus de relations d'\'equivalence de type $\IIi$ ayant la propri\'et\'e T de Kazhdan ont co\^ut 1 --- o\`u le co\^ut d'une relation d'\'equivalence de type $\IIi$ est par d\'efinition le \og 1-covolume \fg\ de cette relation, i.e. l'infimum sur les graphages sym\'etriques $K$ de cette relation du volume $\frac 1 2 \h(K)$ des ar\^etes  de ces graphages (cf \cite{Gaboriau99}).

\bigskip

\textit{Remarque.} Les relations arborables admettent des approximations non triviales, cf. \cite{Gaboriau02,GaboriauPopa04}. Par exemple, \'etant donn\'e un arborage $K$ d'une relation $\R$, toute suite croissante $K_n \subset K$ non essentiellement constante et exhaustant $K$ d\'etermine une approximation non triviale de $\R$. Notons que, la propri\'et\'e T \'etant \'evidemment stable par quotient, une relation d'\'equivalence de Kazhdan ne poss\`ede pas de quotients arborables. (cf. \'egalement \cite{AdamsSpatzier90})

\bigskip
\subsection{Ergodicit\'e forte et familles de Levy.} Nous avons d\'ej\`a constat\'e dans \cite{Pichot04_I} des liens \'etroits entre l'ergodicit\'e forte et le ph\'enom\`ene de concentration de la mesure. Dans ce paragraphe nous caract\'erisons les espaces fortement ergodiques \`a l'aide de familles de Levy qui leurs sont naturellement associ\'ees. 

\bigskip

Notre r\'ef\'erence pour la concentration de la mesure au sens classique est la monographie r\'ecente de Ledoux \cite{Ledoux01}. 

\bigskip

Consid\'erons  d'abord la notion de famille de Levy en termes de fonctions 1-lipschitziennes. Soit $(Y,d,\mu)$ un espace m\'etrique-mesur\'e, i.e.  un espace m\'etrique $(Y,d)$ muni d'une mesure bor\'elienne de probabilit\'e $\mu$. \'Etant donn\'ee une fonction 
\[
f : Y \to \RI
\]
\`a valeurs r\'eelles, on sera plus particuli\`erement int\'eress\'e par les \textit{in\'egalit\'es de concentration} de $f$ autour d'une valeur $m\in \RI$, de la forme
\[
\mu\{\abs{f-m}\leqslant \varepsilon\}\geqslant \delta(\varepsilon)
\]
(o\`u $\delta(\varepsilon)\to 1$ quand $\varepsilon \to \infty$).

\bigskip

\begin{dfn}
Soit $((Y_n,y_n),d_n,\mu_n)_n$ une suite d'espaces m\'etriques-mesur\'es point\'es, o\`u  $y_n\in Y_n$ est le point base de $Y_n$. On suppose que les premiers moments
\[
m_1((Y_n,y_n),d_n,\mu_n)=\int_{Y_n} d_n(y,y_n)d\mu_n(y)\leqslant C<\infty
\]
sont unifom\'ement finis ($C>0$ fix\'e). Nous dirons que  $((Y_n,y_n),d_n,\mu_n)$ forme une \textnormal{famille de Levy} si pour tout $\varepsilon>0$ on a,
\[
\inf_{f} \mu_n\{\abs{f-m}\leqslant \varepsilon\} \to_n 1,
\]  
o\`u l'infimum est pris sur les fonctions 1-lipschitziennes  $f : Y_n\to \RI$  et  $m = \int_{Y_n} fd\mu_n$ est la valeur moyenne de $f$.
\end{dfn}

\bigskip

Dans la suite nous \'etudierons seulement le cas o\`u $(Y_n,d_n) = (H,\norm \cdot)$ est un espace de Hilbert fixe, avec l'origine pour point base. On dira dans ce cas qu'une fonction $f : H \to \RI$ \textit{se concentre au voisinage d'une valeur $m$} lorsque
\[
\mu_n\{\abs{f-m}\leqslant \varepsilon\} \to_n 1
\]

\bigskip

\begin{thm}\label{Levy} Soit $H$ un espace de Hilbert. Soit $\R$ une relation d'\'equivalence ergodique sur un espace de probabilit\'e $(X,\mu)$. Alors $\R$ est fortement ergodique si et seulement si pour toute suite presque invariante $\xi^n : X \to H$ pour la repr\'esentation triviale de $\R$ dans $H$, domin\'ee par une fonction $g\in L^1$, la famille $(H, \norm \cdot, \mu_n)$ est une famille de Levy, o\`u $\mu_n = \xi^n_*\mu$ est la pouss\'ee en avant de $\mu$ sur $H$.
\end{thm}

\begin{dem}
Supposons d'abord que $\R$ soit fortement ergodique et consid\'erons une suite presque invariante domin\'ee 
\[
\xi^n : X \to H.
\]
Soit $\varepsilon >0$ fix\'e. \'Etant donn\'e un entier $n$ consid\'erons une fonction  1-lipschitzienne $f_n : H \to \RI$ telle que 
\[
\mu_n\{f_n-m_n> \varepsilon\}\geqslant \sup_{f} \mu_n\{f-\int fd\mu_n> \varepsilon\} -1/n,
\] 
o\`u $m_n=\int_H f_nd\mu_n$ (et  $\mu_n = \xi^n_*\mu$) et montrons que $\mu_n\{f_n-m_n> \varepsilon\}\to 0$. Supposons qu'il existe une suite extraite de $(f_n)$, encore not\'ee $(f_n)$, telle que 
\[
\mu_n\{f_n-m_n> \varepsilon\}\to \delta>0.
\]
Notons $C^n=\{f_n-m_n> \varepsilon\}\subset H$ et $A_n=\{\xi^n_x \in C^n\} \subset X$. Les champs constants $x\mapsto C^n$ sont bien s\^ur invariant pour la repr\'esentation  triviale. D'apr\`es le lemme \ref{pi-ai} il existe quitte \`a extraire une seconde fois une suite $\alpha=(\alpha_n)$ de nombres r\'eels positifs convergeant vers 0, de sorte que
\[
A_n^\alpha = \{\xi^n_x \in C^n_{\alpha_n}\}
\]
soit asymptotiquement invariante, o\`u $C^n_{\alpha_n}$ est le $\alpha_n$-voisinage de $C^n$ dans $H$. Comme $\R$ est fortement ergodique, cette suite est triviale, i.e.
\[
\mu(A_n^\alpha) \to_n 0~~\mathrm{ou}~~1,
\]
et donc
\[
\mu(A_n^\alpha) \to_n 1.
\]
On a alors
\begin{eqnarray}
m_n=\int_H f_nd\mu_n&=& \int_{A_n^\alpha} f_n\xi^n d\mu+\int_{X\backslash A_n^\alpha}f_n\xi^nd\mu \nonumber \\
&>& \mu(A_n^\alpha)(m_n+\varepsilon-\alpha_n) -\abs{\int_{X\backslash A_n^\alpha} f_n\xi^nd\mu}\nonumber \\
&\geqslant&\mu(A_n^\alpha)(m_n+\varepsilon-\alpha_n)  - \int_{X\backslash A_n^\alpha} g d\mu-\mu(X\backslash A_n^\alpha)\abs{f_n(0)},\nonumber\\
&\geqslant&\mu(A_n^\alpha)(m_n+\varepsilon-\alpha_n)  - \int_{X\backslash A_n^\alpha} g d\mu-\mu(X\backslash A_n^\alpha)(\abs{m_n}+\norm g_1),\nonumber
\end{eqnarray}
o\`u les in\'egalit\'es utilisent le fait que $f_n$ est 1-lipschitz et la d\'efinition de $A_n^\alpha$. Or $\abs{m_n}$ est par d\'efinition born\'e par $\norm g_1$, d'o\`u une contradiction pour $n$ grand. En rempla\c cant $f_n$ par $-f_n$ on a finalement,
\[
\sup_{f} \mu_n\{\abs{f-\int fd\mu_n}> \varepsilon\} \to 0,
\]
d'o\`u le r\'esultat.

\bigskip

R\'eciproquement supposons que $\R$ ne soit pas fortement ergodique et construisons une suite asymptotiquement invariante $\xi^n$ telle que la famille $(H,\norm\cdot,\mu_n)$ associ\'ee ne soit pas une famille de Levy. Soit $(A_n)$ une suite asymptotiquement invariante non triviale pour $\R$.  Soit $(e_n)$ une suite dense de la sph\`ere unit\'e de $H$. (Rappelons que par convention $H$ est  s\'eparable.) Posons
\[
\xi^n = \chi_{A_n}\cdot e_n - \mu(A_n)e_n,
\]
o\`u $\chi_{A_n}$ est la fonction caract\'eristique de $A_n$. On v\'erifie imm\'ediatement que $(\xi^n)$ est une suite presque invariante et domin\'ee (par la fonction constante 1). Soit $\eta\in H$ un vecteur unit\'e.  Consid\'erons la fonction (1-lipschitz) 
\[
\xi \mapsto \abs{\langle \xi\mid\eta\rangle}.
\]
Si  $(H,\norm\cdot,\mu_n)$ \'etait une famille de Levy on aurait en particulier pour tout $\varepsilon>0$,
\[
\mu\{\abs{\langle \xi^n_x\mid\eta\rangle} \leqslant \varepsilon\} \to 1.
\]
Consid\'erons une suite extraite de $(\xi^n)$ (toujours not\'ee $(\xi^n)$) telle que pour la suite $(e_n)$ associ\'ee on ait, 
\[
\abs{\langle e_n\mid \eta\rangle} \to 1.
\]
On a
\[
\abs{\langle \xi^n_x\mid\eta\rangle}=\abs{\chi_{A_n} - \mu(A_n)}\cdot\abs{\langle e_n\mid \eta\rangle}
\]
et donc
\[
\int_X\abs{\langle \xi^n_x\mid\eta\rangle}d\mu = 2\cdot\mu(A_n)(1-\mu(A_n))\abs{\langle e_n\mid \eta\rangle},
\]
d'o\`u une contradiction pour $n$ grand si $\varepsilon$ est suffisamment petit ($\varepsilon\leqslant \delta^2/3$ o\`u $A_n$ est $\delta$-non triviale).
\end{dem}

\bigskip
\subsection{Proximit\'e des champs invariants et presque invariants.}

L'un des corollaires importants des r\'esultats techniques du paragraphe \ref{pinv} est le th\'eor\`eme  suivant,  dont l'analogue p\'eriodique est bien connu, cf. \cite{HarpeValette89}.  La d\'emonstration de ce th\'eor\`eme comporte 3 \'etapes (th\'eor\`eme \ref{prox}, lemme \ref{triv}, lemme \ref{trivinf}).

\bigskip

\begin{thm}\label{prox}
Soit $\R$ une relation d'\'equivalence ergodique de type $\IIi$ ayant la propri\'et\'e T. Soit $\pi$ une repr\'esentation  de $\R$ poss\'edant une suite $\xi^n$ presque invariante telle que $\norm{\xi^n_x}_x\leqslant g(x)$ pour une fonction $g\in L^1(X)$. Il existe quitte \`a extraire une suite $\zeta_n$ de champs invariants domin\'es par $g$ tels que 
\[
\norm {\xi_x^n -  \zeta^n_x}_x \to 0
\]
pour presque tout $x\in X$.
\end{thm}

\begin{dem} Soit $\varepsilon>0$. Soit $\pi$ une repr\'esentation de $\R$ sur un champ hilbertien $H$ de base $X$. D\'ecomposons $H = H^t+\tilde H$ en somme orthogonale, o\`u $H^t$ est le champ mesurable et stable par $\pi$ engendr\'e par les champs invariants. Consid\'erons la famille  
\[
P_x : H_x \to H_x
\]
$x\in X$, de projecteurs orthogonaux sur $\tilde H$. Pout tout $(x,y)\in \R$ on a
\[
P_x\pi(x,y)=\pi(x,y)P_y
\]
Soit $\xi^n$ une suite presque invariante de champs de vecteurs tels que $\norm{\xi^n_x}_x \leqslant g(x)$ avec $g\in L^1$, et $\xi^n = \zeta^n+\tilde \xi^n$ la d\'ecomposition de $\xi^n$ dans $H =H^t+\tilde H$. On a
\[
\norm{\pi(x,y)\tilde \xi^n_y-\tilde \xi^n_x}= \norm{P_x(\pi(x,y)\xi^n_y-\xi^n_x)}\leqslant\norm{\pi(x,y)\xi^n_y-\xi^n_x},
\]
donc $\tilde \xi^n$ (et de m\^eme $\zeta^n$) est une suite presque invariante. Par construction, $\pi_{|\tilde H}$ ne contenant pas de champs presque invariants non triviaux, et puisque  $\norm {\tilde \xi^n_x}_x \leqslant g(x)$, on a, d'apr\`es le th\'eor\`eme \ref{dom}, $\norm{\tilde \xi^n}_1 \to_n 0$, i.e. $\norm{\tilde \xi^n_x}_x \to_n 0$ presque s\^urement quitte \`a extraire. 
Par suite 
\[
\norm {\xi_x^n -  \zeta^n_x}_x \to 0
\]
pour presque tout $x\in X$.  Soit $T$ un champ d'isomorphismes entrela\c cant la restriction $\pi_{|H^t}$ de $\pi$ \`a $H^t$ et la repr\'esentation triviale de $\R$, agissant sur un champ constant $X\times H^\1$ de base $X$. Consid\'erons le champ d\'efini par
\[
\tilde \zeta^n_x = T_x(\zeta_x^n) \in  H^\1.
\]
Alors $\tilde \zeta^n : X\to H^\1$ est une suite presque invariante domin\'ee. 
Soit $\eta^n$ le champ d\'efini par
\[
\eta_x^n = \tilde \zeta_x^n -\int_X\tilde \zeta_x^nd\mu.
\]
Le lemme suivant (lemme \ref{triv}) montre que $\eta_x^n \to 0$ presque s\^urement quitte \`a extraire. 
Il en r\'esulte que 
\[
\norm {\xi_x^n - T_x^{-1}\int_X\tilde \zeta_x^nd\mu}_x\leqslant \norm {\xi_x^n -  \zeta^n_x}_x + \norm{\tilde \zeta_x^n -\int_X\tilde \zeta_x^nd\mu}_x \to 0
\]
pour presque tout $x\in X$, d'o\`u le r\'esultat. 
\end{dem}

\bigskip

\begin{lem}\label{triv}
Soit $\R$ une relation d'\'equivalence fortement ergodique sur un espace de probabilit\'e $(X,\mu)$. On consid\`ere une suite $(\xi^n : X \to H)_n$ presque invariante pour la repr\'esentation triviale de $\R$ sur un espace de Hilbert $H$, domin\'ee par une fonction $g\in L^1(X)$, et telle que
\[
\norm{\int_X \xi^nd\mu}\to_n 0.
\]
Si $H$ est de dimension finie, ou si $\R$ est de type fini et pr\'eserve la mesure $\mu$, alors il existe une suite extraite de  $(\xi^n)$ qui converge presque surement vers 0. 
\end{lem}

\begin{dem} Si $H$ est de dimension fini, il s'agit d'un corollaire imm\'ediat du th\'eor\`eme \ref{Levy} de concentration. En effet choisissons pour fonctions lipschitziennes les fonctions coordonn\'ees de $H$ (consid\'er\'e comme espace de Hilbert r\'eel), qui se concentrent au voisinage de 0 du fait que leurs valeurs moyennes tend vers 0. Plus pr\'ecis\'ement on a
\[
\inf_\eta \mu\{\abs{\langle \xi^n_x\mid\eta\rangle}\leqslant \varepsilon\}\to_n 1.
\]
o\`u $\eta$ parcourt les vecteurs unit\'es de $H$. Par suite \'etant donn\'e $\varepsilon>0$, on en d\'eduit (si $H$ est de dimension finie) que
\[
\mu(O_\varepsilon^n)\to_n 1,
\]
o\`u $O_\varepsilon^n$ est l'ensemble des $x$ tels que $\xi^n_x$ est dans la boule de centre 0 et de rayon $\varepsilon$ de $H$,
\begin{eqnarray}
\norm{\xi^n} &=& \int_{O_\varepsilon^n} \norm{\xi_x^n} d\mu(x) + \int_{X\backslash O_\varepsilon^n} \norm{\xi_x^n} d\mu(x)\nonumber\\
&\leqslant& \varepsilon + \int_{X\backslash O_\varepsilon^n} g d\mu\nonumber
\end{eqnarray}
donc $\norm{\xi^n}\leqslant 2\varepsilon$ pour $n$ grand (il en r\'esulte quitte \`a extraire que $\xi^n$ converge vers 0 presque s\^urement).  Suivant la comparaison de M. Gromov \cite[page 141]{Gromov00_SQ}, le ph\'enom\`ene concentration concerne \textit{a priori} principalement les \og observables\fg\ $f:H \to Y$ o\`u $Y$ est un \og \'ecran\fg\ de basse dimension. Pour \'etudier le cas o\`u $Y$ est (hilbertien) de dimension infinie, nous utiliserons des arguments spectraux (lemme \ref{trivinf}), ce qui permettra de conclure la preuve de ce lemme.
\end{dem}

\bigskip

Nous avons donc montr\'e, modulo le lemme  \ref{trivinf}, le th\'eor\`eme suivant.

\begin{thm}\label{caract-prox}
Soit $\R$ une relation d'\'equivalence ergodique de type $\IIi$ sur un espace de probabilit\'e $(X,\mu)$. Alors $R$ poss\`ede la propri\'et\'e T de Kazhdan si et seulement si pour toute repr\'esentation $\pi$ de $\R$ sur un champ hilbertien $H$ de base $X$ et toute suite presque invariante $\xi^n$ de $H$ domin\'e par une fonction $g\in L^1(X)$, il existe quitte \`a extraire une suite $\zeta^n$ de champs invariants domin\'es par $g$ tels que 
\[
\norm {\xi_x^n -  \zeta^n_x}_x \to 0
\]
pour presque tout $x\in X$. 
\end{thm}

\bigskip
\subsection{Groupes discrets et relations d'\'equivalence.} Nous montrons dans ce paragraphe qu'un groupe discret agissant librement en pr\'eservant une mesure de probabilit\'e poss\`ede la propri\'et\'e T si et seulement si la relation d'\'equivalence engendr\'ee la poss\`ede. Ce r\'esultat a \'et\'e obtenu par Zimmer, dans le cas des actions faiblement m\'elangeantes \cite{Zimmer84}, puis  par Popa  (par des techniques d'alg\`ebres de von Neumann \cite{Popa86}) et plus r\'ecemment par Anantharaman-Delaroche  (par des techniques cohomologiques \cite{Anan04}) dans le cas des actions ergodiques.

\bigskip

\begin{thm}
Soit $\R=\R_\alpha$ une relation obtenue par action d'un groupe discret $\G$ pr\'eservant une mesure de probabilit\'e. Si $\G$ a la propri\'et\'e T, alors $\R$ l'a \'egalement. R\'eciproquement si $\R$ a la propri\'et\'e T, et si l'action $\alpha$ est essentiellement libre et ergodique, alors $\G$ l'a \'egalement.
\end{thm}

\begin{dem}
Soit $\G$ un groupe de Kazhdan, $\R=\R_\alpha$ une relation obtenue par action de $\G$ pr\'esevant une mesure de probabilit\'e $\mu$, et $\pi$ une repr\'esentation de $\R$ poss\'edant presque des champs invariants ; alors la repr\'esentation $\overline \pi$ de $\G$ sur $L^2(X,H)$ obtenue en int\'egrant $\pi$ poss\`ede presque des vecteurs invariants. En effet, soit $F \subset \G$   une partie finie, $\varepsilon >0$,  $K = \{\graph(\gamma^{-1})\}_{\gamma \in F}\subset \R$, et une suite $\xi^n \in L^2$ de champs $(\varepsilon_n, K_n)$-invariant (o\`u $\varepsilon_n \to 0$ et $K_n$ est exhaustive croissante) tels que $\norm {\xi^n} _{L^2} =1$ et $\norm{\xi^n_x}\leqslant g(x)$ pour une fonction positive $g\in L^2$. Soit un entier $N$ tel que
\[
 \int_{K\backslash K_N} (g(x) + g(y))^2 d\h(x,y) \leqslant \varepsilon^2/2
\]
et $\varepsilon_N^2 \leqslant \varepsilon^2/2$. On a pour tout $\gamma \in F$
\begin{eqnarray}
\norm {\overline\pi(\gamma) \xi^N - \xi^N}_2^2 &=& \int_X \norm {\pi(x,\gamma^{-1} x)\xi^N_{\gamma^{-1} x} - \xi^N_x}^2 d\mu(x) \nonumber \\
&\leqslant& \varepsilon^2/2 + \int_{\graph(\gamma^{-1}) \cap K_N}\norm { \pi(x,y)\xi^N(y) - \xi^N(x)}^2 d\h(x,y) \nonumber \\
&\leqslant& \varepsilon^2/2 + \varepsilon_N^2 \leqslant \varepsilon^2. \nonumber
\end{eqnarray}
Il existe donc une fonction non triviale $\xi \in L^2(X,H)$ invariante par $\G$ (i.e. v\'erifiant $\pi(x,\gamma^{-1}x)\xi_{\gamma^{-1}x}=\xi_x$ pour presque tout $x\in X$). Comme $\G$ engendre $\R$, cette fonction est invariante pour $\pi$ et  $\R$ a la propri\'et\'e T.

\medskip

R\'eciproquement soit $\pi : \G \to U(H)$ une repr\'esentation de $\G$ dans $H$ ayant presque des vecteurs invariants. Soit $\tilde \pi$ la repr\'esentation  de $\R$ sur le champ constant d'espaces $H_x =H$ et de base $X$, d\'efinie par l'expression
\[
(x,y) \mapsto (h\in H_y \mapsto \pi(\gamma^{-1})h\in H_x),
\]
o\`u  $\gamma$ est l'unique \'el\'ement de $\Gamma$ v\'erifiant l'\'equation $\alpha(\gamma)(x)=y$ ($\alpha$ \'etant suppos\'ee libre). Tout vecteur $(F,\varepsilon)$-invariant pour $\pi$, disons $\xi$, d\'efinit un champ constant $(K,\varepsilon)$-invariant pour $\tilde \pi$,
\[
x\mapsto \xi
\]
o\`u $K = \{\graph(\gamma^{-1})\}_{\gamma \in F}$. Consid\'erons une suite presque invariante pour $\pi$ et notons $\xi^n : X \to H$ la suite de champs constants associ\'es. $\R$ ayant la propri\'et\'e T, il existe une suite  $(\zeta^n)$  de champs invariants domin\'es tels que $\zeta^n_x -\xi^n_x \to 0$ presque s\^urement (cf. th. \ref{prox}). Soit
\[
\overline\zeta_n = \int_X \zeta^n_x d\mu(x)=\sum_i \int_X   \langle \zeta^n_x\mid e_i\rangle d\mu(x)e_i
\]
la valeur moyenne de $\zeta^n$, o\`u $(e_i)$ est une base hilbertienne de $H$. Alors
\[
\pi(\gamma)\overline\zeta_n = \pi(\gamma) \int_X \zeta^n d\mu = \int_X \pi(\gamma) \zeta^n_x d\mu(x)=\int_X \zeta^n_{\alpha(\gamma)x} d\mu(x)=\overline \zeta_n
\]
(la deuxi\`eme in\'egalit\'e r\'esulte de la continuit\'e de $\pi(\gamma)$ et la derni\`ere d'un changement de variable, $\mu$ \'etant invariante). Or
\[
\norm {\overline \zeta_n -\xi^n} = \norm{\int_X\zeta^n-\xi^nd\mu} \leqslant \int_X\norm{\zeta^n-\xi^n}d\mu \to 0
\]
donc $\overline\zeta_n$ est non trivial pour $n$ grand. Ainsi $\G$ a la propri\'et\'e T.
\end{dem}

\bigskip

Le r\'esultat suivant a \'et\'e d\'emontr\'e par Furman \cite{Furman99}. Il r\'esulte aussi imm\'ediatement du th\'eor\`eme ci-dessus et de l'invariance par isomorphisme stable de la propri\'et\'e T pour les relations d'\'equivalence.

\medskip

\begin{cor}
La propri\'et\'e T est un invariant d'\'equivalence mesurable de groupes discrets.
\end{cor}

Rappelons que deux groupes discrets $\G$ et $\Lambda$ sont dit \textit{mesurablement \'equivalents} s'il existe un espace bor\'elien standard $\Omega$ muni d'une mesure $\sigma$-finie $\h$, et des actions de $\G$ et $\Lambda$ sur $\Omega$ qui soient libres et commutantes, qui pr\'eservent la mesure $\h$, et qui admettent chacune un domaine fondamental de mesure finie (cette d\'efinition est d\^ue \`a M. Gromov). Par exemple, deux r\'eseaux de covolume fini d'un m\^eme groupe de Lie sont mesurablement \'equivalent (agissant par multiplication \`a gauche et \`a droite).

Deux groupes discrets sont mesurablement \'equivalents si et seulement s'ils admettent deux actions libres de type $\IIi$ stablement isomorphes (\cite[\textsection 2]{Furman99} et \cite[\textsection 6]{Gaboriau02}).

\vspace{1cm}
\section{Diffusion hilbertienne et in\'egalit\'es de Poincar\'e}\label{Diffus}

Soit $H$ un espace de Hilbert.\\

\'Etant donn\'e un op\'erateur hermitien born\'e $D$ sur $H$ de norme $\leqslant 1$ (contraction), on consid\`ere l'op\'erateur positif
\[
\Delta_p = \id - D^p
\]
o\`u $\id$ est l'identit\'e et $p \geqslant 1$. On note
\[
E_p(\xi) = E_{D,p}(\xi) = \langle \Delta_p \xi\mid \xi \rangle
\]
l'energie (de Dirichlet) de $\xi\in H$ relative \`a la diffusion $\xi \mapsto D^p\xi$. 

\bigskip

Les points fixes de la diffusion $D$, i.e. les vecteurs $\xi$ v\'erifiant $D\xi = \xi$, sont les vecteurs dont l'\'energie $E=E_1$ est nulle.\\

\textit{Exemple (Marche al\'eatoire sur la sph\`ere hilbertienne $\SI^\infty \subset H$).} Soit $\G$ un groupe de type fini. On fixe une marche al\'eatoire invariante $\nu$ sur $\G$, de support un syst\`eme g\'en\'erateur fini \textit{sym\'etrique} $S$ de $\G$. Ainsi $\nu$ consiste en la donn\'ee de $\#S$ nombres r\'eels strictement positifs
\[
\nu(e\to s)
\]
de somme 1, qui d\'eterminent par invariance la probabilit\'e $\nu(\gamma\to s\gamma)$ d'aller de $\gamma$ \`a $s\gamma$ en 1 pas. On suppose que $\nu$ est sym\'etrique, au sens o\`u $\nu(e\to s) = \nu(e\to s^{-1})$.
 Soit $\pi$ une repr\'esentation unitaire  de $\G$ sur un espace de Hilbert $H$. L'op\'erateur
\[
D_{\nu,\pi} = \sum_S \nu(e\to s) \pi(s)
\]
de diffusion associ\'e \`a la marche al\'eatoire sur la sph\`ere unit\'e de $H$ est alors hermitien de norme $\leqslant 1$, sans point fixe si et seulement si $\pi$ est sans point fixe (en effet le barycentre $D_{\nu,\pi}(\xi)$ des vecteurs unitaires pond\'er\'es $(\pi(s)\xi,\nu(e\to s))$ est de norme 1 si et seulement si $\pi(s)\xi = \xi$ pour tout $s$). Partant d'un vecteur unitaire $\xi$ de $H$, on se d\'eplace  en $\pi(s)\xi$ avec probabilit\'e $\nu(e\to s)$. On note
\[
\nu^2(\gamma\to \gamma')=\nu*\nu(\gamma\to \gamma') = \sum_{\tau\in \G} \nu(\gamma \to \tau)\nu(\tau \to \gamma')
\]
la probabilit\'e d'aller de $\gamma$ \`a $\gamma'$ en 2 pas sur $\G$, et de m\^eme $\nu^n= \nu^{n-1}*\nu$. On a $D_{\nu^n,\pi} = D_{\nu,\pi}^n$ pour les diffusions hilbertiennes associ\'ees \`a $\nu$.

\bigskip
\bigskip

Soit $D$ une contraction de $H$. Rappelons que $\Sp(D)\subset [-1,1]$, o\`u $\Sp(D)$ est le spectre de $D$. On note 
\[
\kappa = \kappa_D=\sup \{y \in \Sp(D), y \neq 1\}.
\]
La pr\'esence d'un \og trou\fg\ 
\[
\lambda = 1 - \kappa >0
\]
dans $\Sp(D)$ \'equivaut \`a la pr\'esence  d'un trou (de m\^eme taille) dans le spectre des \'energies, i.e. 
\[
E(\xi) = E_1(\xi) \geqslant \lambda >0,
\]
pour tout $\xi \in \SI^\infty$ non fixe.

\bigskip
\bigskip

{\bf In\'egalit\'e de Dirichlet.} \textit{Soit $D$ une contraction. Alors $\kappa <1$ si et seulement s'il existe une constante $c_\infty < \infty$ telle que
\[
\norm{\xi-\overline \xi}^2 \leqslant c_\infty E(\xi)
\]
pour tout $\xi \in H$ (o\`u $\overline \xi$ est la projection orthogonale de $\xi$ sur les points fixes de $D$). La valeur optimale de cette constante est $c_\infty = 1/(1-\kappa)=1/\lambda$.}

\bigskip

\begin{dfn}
Les constantes  $c_\infty$ et $\lambda$  sont appel\'es \textnormal{constantes de relaxation} de la diffusion $D$.
\end{dfn}

\bigskip

\bigskip

{\bf In\'egalit\'es de Poincar\'e.} \textit{Soit $D$ une contraction. Alors $\kappa <1$ si et seulement s'il existe $n\geqslant 2$ et une constante $c_n < n$ tels que
\[
E_n(\xi) \leqslant c_n E(\xi)
\]
pour tout $\xi \in H$. L'in\'egalit\'e $c_n<n$ est alors vraie pour tout $n\geqslant 2$. La valeur optimale de la constante $c_n$ est $c_n = 1 + \kappa+ \kappa^2 + \ldots + \kappa^{n-1}$.}\\

\begin{dem} On peut supposer, quitte \`a consid\'erer l'orthogonal des points fixes, que $D$ est sans point fixe. L'in\'egalit\'e revient \`a dire que l'op\'erateur $\Delta T$ est positif, o\`u
\[
\Delta= \Delta_1 = \id- D \ \mathrm{et} \  T = c_n - (\id + D + \ldots +D^{n-1}).
\] 
Or ceci est \'equivalent \`a la positivit\'e de $T$. En effet $\Delta$ est positif et injectif donc d'image dense, et en \'ecrivant $\xi = \lim \sqrt {\Delta} \xi_n$, on obtient
\[
\langle T \xi \mid \xi \rangle = \lim_n\, \langle T\sqrt {\Delta}\xi_n\mid \sqrt {\Delta}\xi_n\rangle = \lim_n \langle \Delta T  \xi_n \mid \xi_n \rangle \geqslant 0.
\]
La r\'eciproque est \'evidente. On a pour $0\leqslant \kappa\leqslant 1$,
\[
D \leqslant \kappa \ssi \id + D + \ldots +D^{n-1}\leqslant 1 + \kappa + \kappa^2 + \ldots + \kappa^{n-1}
\]
et $\kappa <1 \ssi 1+ \kappa + \kappa^2 + \ldots + \kappa^{n-1}<n$.
\end{dem}

\bigskip

\begin{dfn}
Les constantes  $c_n$ appel\'ees \textnormal{constantes de Poincar\'e} de la diffusion $D$.
\end{dfn}

\bigskip
\bigskip

Soit $\G$ un groupe de type fini.\\

Fixons comme dans l'exemple ci-dessus un syst\`eme g\'en\'erateur $S$ sym\'etrique fini de $\G$ et $\nu$ une marche al\'eatoire invariante sym\'etrique de support $S$. \'Etant donn\'ee une repr\'esentation $\pi$ de $\G$ sur un espace de Hilbert $H$, on note $D_{\nu,\pi}$ la marche al\'eatoire associ\'ee sur $\SI^\infty \subset H$, et $\kappa(\nu,\pi)$, $c_n(\nu,\pi)$ les constantes correspondantes.\\

Notons que (par des consid\'erations barycentriques \'evidentes) $\G$ a la propri\'et\'e T si et seulement si pour toute repr\'esentation $\pi$,  on a
\[
\kappa(\nu,\pi) <1.
\]
Les in\'egalit\'es de Poincar\'e et Dirichlet s'\'ecrivent
\[
c_n(\nu,\pi) < n.
\]
On voit alors facilement du fait qu'elles sont atteintes que les constantes $\kappa(\nu,\pi)$ et $c_n(\nu,\pi)$ sont alors uniformes en $\pi$ (consid\'erer une suite $\pi_n$ de repr\'esentations pour lesquelles la constante tend vers la valeur maximale et faire la somme directe de ces repr\'esentations). En d'autres termes si l'on note
\[
\kappa(\G,\nu)= \sup_\pi \kappa(\nu,\pi) \quad \mathrm{et} \quad c_n(\G,\nu) = \sup_\pi c_n(\nu,\pi),
\]
alors \textit{$\G$ poss\`ede la propri\'et\'e T de Kazhdan si et seulement si l'une des in\'egalit\'es \'equivalentes $\kappa(\G,\nu) < 1$ ou $c_n(\G,\nu) < n$ est v\'erifi\'ee ($n\geqslant 2$).}

\bigskip

Nous renvoyons \`a \cite{Gromov03,Ghys03} pour ce qui pr\'ec\`ede.

\vspace{1cm}
\section{Le point de vue spectral}\label{dem}

Le but de ce paragraphe est la d\'emonstration des th\'eor\`emes spectraux \'enonc\'es au cours de l'introduction (th. \ref{ef-intr} et \ref{T-intr}).

\bigskip

\subsection{Marche al\'eatoire sur une relation d'\'equivalence.} Soit $\R$ une relation d'\'equi-valence mesur\'ee sur un espace de probabilit\'e $(X,\mu)$.\\

Une \textit{marche al\'eatoire sur (les orbites de) $\R$} est la donn\'ee d'une famille mesurable de mesures de probabilit\'e $(\nu_x)_{x\in X}$ telle que $\nu_x$ soit support\'ee sur l'orbite de $x$. Plus pr\'ecis\'ement, $\nu : \R \to [0,1]$ est une fonction mesurable d\'efinie sur $\R$ telle que la somme de chaque fibre horizontale $\R^x$ soit 1. On notera $\nu(x\to y)=\nu(x,y)$. (Observons que $\nu$ s'\'etend naturellement par \'equivariance en une marche al\'eatoire $\tilde \nu$ sur l'ensemble d\'enombrable quasi-p\'eriodique $\R$, au sens de \cite{Pichot04_I}, en posant $\tilde \nu((x,y)\to(x,z))=\nu(y\to z)$. Nous avons choisi ici de travailler avec $\nu$ plut\^ot que $\tilde \nu$  pour simplifier les notations.)

\bigskip

\begin{dfn} Une marche al\'eatoire $\nu$ est dite \textnormal{sym\'etrique relativement \`a $\mu$} si
\[
\nu(x\to y)\sqrt{\delta(x,y)} =\nu(y\to x)\sqrt{\delta(y,x)}.
\]
o\`u $\delta$ est la fonction modulaire de $\mu$.
\end{dfn}

\bigskip

Rappelons que $\delta$ v\'erifie l'\'equation
\[
d\h(x,y) = \delta(x,y) d\h^{-1}(x,y),
\]
o\`u $(y,x)=(x,y)^{-1}$ est l'inversion (ainsi $\h^{-1}$ est la mesure de d\'ecompte vertical).

\bigskip

\textit{Exemple.} Soit $\R$ une relation d'\'equivalence mesur\'ee sur un espace bor\'elien standard $(X,\mu)$. Soit $K$ un graphage sym\'etrique u.l.b. de $\R$. La \textit{marche al\'eatoire r\'eguli\`ere sur $K$} est d\'efinie par
\[
\nu_K(x\to y) = \frac{\sqrt{\delta(y,x)}}{\sum_{(x,y)\in K^x} \sqrt{\delta(y,x)}}
\]
si $(x,y) \in K$ et 0 sinon. Elle est sym\'etrique relativement \`a la mesure $\tilde \mu$ d\'efinie par
\[
d\tilde \mu (x) = \delta(x) d\mu(x)
\]
(\'equivalente \`a $\mu$) o\`u $\delta(x)= \sum_{(x,y)\in K^x} \sqrt{\delta(y,x)}$. En effet on a \[
\tilde\delta(x,y)= \delta(x,y)\frac{\delta(x)}{\delta(y)}.
\]

\bigskip

\begin{dfn}
On dit qu'une marche al\'eatoire $\nu$ sur les orbites de $\R$ est \textnormal{sym\'etrique born\'ee} si $\nu$ est sym\'etrique relativement \`a $\mu$, si son support 
\[
K=\mathrm{supp}(\nu)=\{ (x,y)\in\R\mid \nu(x\to y) >0\}
\]
est un \textnormal{graphage} sym\'etrique de $\R$, et s'il existe un nombre r\'eel $\eta > 0$ v\'erifiant $\nu(x\to y) \geqslant \eta$ pour presque tout $(x,y) \in K$.
\end{dfn}

\medskip

Pour tout graphage sym\'etrique u.l.b. $K$ de $\R$, la marche al\'eatoire r\'eguli\`ere $\nu_K$ associ\'ee est sym\'etrique born\'ee.\\

\bigskip

\textit{Int\'egration d'une marche al\'eatoire \`a coefficients dans une repr\'esentation.} Fixons une marche al\'eatoire $\nu$ sur les orbites de $\R$ sym\'etrique relativement \`a $\mu$, de support un graphage sym\'etrique $K$ de $\R$. \'Etant donn\'ee  une repr\'esentation $\pi$ de $\R$ sur un champ d'espaces de Hilbert $H$ de base $X$, on construit une famille mesurable de marches al\'eatoires op\'erant sur les sections mesurables de $H$ en posant
\[
(D_{\nu,\pi}\xi)_x = \sum_{y\sim x} \nu(x\to y)\pi(x,y)\xi_y.
\]
Cette famille d\'efinit un op\'erateur hermitien $D_{\nu,\pi}$ sur l'espace $L^2(X,H)$ des sections de carr\'e int\'egrable sur $X$ ; en effet
\[
\langle D_{\nu,\pi}\xi \mid \eta \rangle = \int_K \langle \nu(x\to y)\pi(x,y)\xi_y\mid \eta_x \rangle d\h(x,y) = \langle \xi\mid D_{\nu,\pi}\eta\rangle.
\]
De plus $D_{\nu,\pi}$ est born\'e de norme $\leqslant 1$ du fait que, pour presque tout $x$, 
\[
\norm{(D_{\nu,\pi}\xi)_x}_x^2 \leqslant \sum_{y\sim x} \nu(x\to y)\norm{\xi_y}_y^2.
\]
On obtient ainsi une diffusion hilbertienne, dont on note $\kappa_\pi$ la plus grande valeur spectrale non triviale (cf \textsection\ref{Diffus}).\\

\begin{dfn}
On dira que $D_{\nu,\pi}$ est la \textnormal{diffusion hilbertienne} associ\'ee \`a $\nu$ et $\pi$.
\end{dfn}

\bigskip

\textit{Constantes de Poincar\'e.} Si alors $\xi \in L^2(X,H)$ est un champ de carr\'e int\'egrable sur $X$ on note 
\[
E_{\pi}(\xi) = E_{D_{\nu,\pi}}(\xi) = \norm{\xi}^2 - \langle D_{\nu,\pi}\xi \mid \xi \rangle
\]
et
\[
E_{\pi,2}(\xi) = E_{D_{\nu,\pi^2}}(\xi) = \norm{\xi}^2 - \norm{ D_{\nu,\pi}\xi }^2,
\]
et on appelle (seconde) constante de Poincar\'e de $\pi$ associ\'ee \`a $\nu$ la plus petite constante $c_2(\pi)$ v\'erifiant
\[
E_{\pi,2}(\xi) \leqslant c_2(\pi) E_{\pi}(\xi).
\]

\vspace{.4cm}

\textit{Remarque terminologique.} Une \textit{diffusion} sur un espace $X$ est, au sens de \cite{Gromov03}, une application $x\mapsto \nu_x=\nu(x\to \cdot)$ de $X$ vers les mesures de probabilit\'e sur $X$, e.g. une marche al\'eatoire. Une \textit{codiffusion} sur un espace $H$ est une application $c$ des mesures de probabilit\'e sur $H$ vers $H$, telle que l'image d'une mesure de Dirac $\delta_\xi$ soit $\xi$ et telle que $c^{-1}(\xi)$ soit convexe pour tout $\xi$ (cf. \cite{Gromov03}). Lorsque $H$ est un espace de Hilbert, la codiffusion naturelle est l'application (affine) de centre de masse
\[
c(\nu) = \int_{H} \xi d\nu(\xi).
\]
\'Etant donn\'ee une repr\'esentation de $\R$ sur un champ hilbertien $H$, on obtient un op\'erateur lin\'eaire qui \`a un champ de vecteurs $\xi$ associe 
\[
x \mapsto c_x((\overline \xi_x)_* (\nu_x))
\]
o\`u $\overline \xi$ est l'extension \'equivariante de $\xi$ \`a (l'ensemble d\'enombrable quasi-p\'eriodique) $\R$ d\'efinie par $\overline \xi_x(y) = \pi(x,y)\xi_y$ pour $x\sim y$, et $c_x$ est la codiffusion naturelle sur $H_x$. C'est cet op\'erateur, restreint aux champs de carr\'e int\'egrable, que nous avons appel\'e diffusion.

\bigskip
\bigskip

\textit{Gradient d'un champ de vecteurs.} \'Etant donn\'e un champ de vecteurs
\[
\xi : X \to H
\]
on d\'efinit son gradient 
\[
d\xi : \R \to H
\]
par l'expression
\[
d\xi(x,y) = \pi(x,y)\xi_y -\xi_x.
\]

\bigskip

Une marche al\'eatoire $\nu$ sur $\R$ d\'efinit canoniquement une mesure de probabilit\'e $\h_\nu$ sur $\R$ par l'expression
\[
\h_\nu(K) =\int_X \sum_{y\in K^x}\nu(x\to y) d\mu(x)
\]
(o\`u $K$ est une partie bor\'elienne de $\R$). Si $f : \R \to H$ est une fonction de carr\'e int\'egrable pour $\h_\nu$ on pose
\[
\norm f _\nu^2 = \int_\R \norm {f(x,y)} ^2 d\h_\nu(x,y) = \int_X \sum_{(x,y)\in \R} \norm {f(x,y)} ^2 \nu(x\to y) d\mu(x).
\]
Alors pour tout champ $\xi \in L^2(X,H)$ de carr\'e int\'egrable, $d\xi \in L^2(\R,H,\h_\nu)$ et on a
\[
E_\pi(\xi) = \frac 1 2 \norm{d\xi}^2_\nu = \frac 1 2 \int_X \sum_{(x,y)\in \R} \norm {\pi(x,y)\xi_y-\xi_x} ^2 \nu(x\to y) d\mu(x)
\]
qui coincide donc avec l'\'energie locale moyenne 
\[
E_\pi(\xi)=\int_X E_\pi(\xi,x) d\mu(x)
\]
o\`u
\[
E_\pi(\xi,x)=\frac 1 2 \sum_{(x,y)\in \R} \norm {\pi(x,y)\xi_y-\xi_x} ^2 \nu(x\to y).
\]

\bigskip

\subsection{Caract\'erisations spectrales.} Commen\c cons par montrer le lemme suivant. 

\bigskip

\begin{lem} \label{CFW} Soit $\R$ une relation d'\'equivalence ergodique. On fixe une marche al\'eatoire $\nu$ sym\'etrique relativement \`a une mesure de probabilit\'e $\mu$, dont le support $K$ engendre $\R$.  Soit $\pi$ une repr\'esentation de $\R$.

i. $\pi$ contient des champs invariants non triviaux si et seulement si $D_{\nu,\pi}$ a des points fixes non triviaux.

ii. Si $\kappa_\pi < 1$, alors  pour toute suite presque invariante $\xi^n$ telle que $\norm{\xi^n_x}_x\leqslant g(x)$ pour une fonction $g\in L^2(X)$, il existe (\`a extraction pr\`es) une suite $\zeta^n$ de champs invariants tels que 
\[
\norm {\xi_x^n -  \zeta^n_x}_x \to 0
\]
pour presque tout $x\in X$.

iii. R\'eciproquement, si $\nu$ est born\'ee, si $K$ ne contient pas de suites de F\o lner \'evanescentes (cf. \cite{Pichot04_I}), et si pour toute suite presque invariante $\xi^n$ telle que $\norm{\xi^n}_2=1$ et $\norm{\xi^n_x}_x\leqslant g(x)$ pour une fonction $g\in L^2(X)$, il existe (\`a extraction pr\`es) une suite $\zeta^n$ de champs invariants tels que 
\[
\norm {\xi_x^n -  \zeta^n_x}_x \to 0
\]
pour presque tout $x\in X$, alors $\kappa_\pi <1$.
\end{lem}

\medskip

\begin{dem}
\textit{i.} Soit $\xi$ un champ fixe par $D_{\nu,\pi}$. Alors $E_\pi(\xi) = 0$, donc $\pi(x,y)\xi_y=\xi_x$ pour presque tout $(x,y)$ dans $K=\mathrm{supp}(\nu)$. Comme $K$ engendre $\R$, $\xi$ est un champ invariant.\\

\textit{ii.}  Montrons que si la conclusion est fausse, alors $\kappa_\pi = 1$.
Notons $V\subset L^2(X,H)$ l'orthogonal dans $L^2(X,H)$ du sous-espace engendr\'e par les vecteurs invariants. Par hypoth\`ese il existe, quitte \`a extraire, une suite $\xi_n \in V$ de champs de vecteurs $(\varepsilon_n, F_n)$-invariants (o\`u $\varepsilon_n \to 0$ et $\h_1(F_n) \to 0$) qui soit uniform\'ement born\'ee par une fonction $g\in L^2$ et telle que $\norm{\xi_n}_2 =1$. Alors $\h_\nu(F_n) \to 1$ et
\[
E_\pi(\xi_n) \leqslant \varepsilon_n \int_{F_n} \nu(x\to y) d\h(x,y) +  \int_{K\backslash F_n} (g(x) + g(y))^2d\h_\nu(x,y).
\]
Donc $E_\pi(\xi_n) \to 0$. Ceci entra\^\i ne que $\kappa_\pi =1$.\\

\textit{iii.} R\'eciproquement supposons par l'absurde que $\kappa_\pi = 1$ et obtenons une contradiction. Consid\'erons  donc une suite $\xi^n\in V$ de champs de vecteurs $X \to H$, de norme $\norm{\xi^n}_2 =1$, et dont l'\'energie tend vers 0. Quitte \`a extraire on peut supposer cette suite presque invariante. En effet pour tout $k$ l'\'energie de $\xi^n$ relativement \`a la marche al\'eatoire en $k$ pas (de support $K^k$) converge vers 0 avec $n$, et cette marche \'etant born\'ee pour tout $k$, on conclut par extraction diagonale.
  
Fixons $\eta \geqslant 1$. D'apr\`es le lemme \ref{pi-ai}, il existe quitte \`a extraire une suite $(\eta_n)$ de nombres r\'eels convergant vers $\eta$ (en d\'ecroissant), telle que la suite $(A_n)$ d\'efinie par
\[
A_n = \{ \frac 1 \eta_n \leqslant \norm {\xi^n} \leqslant \eta_n\}.
\]
soit asymptotiquement invariante, et de mesure convergeant vers $\delta \in [0,1]$.  Comme $K$ ne contient pas de suites de F\o lner \'evanescentes, $\R$ est fortement ergodique, donc cette suite est triviale, i.e. $\delta = 0$ ou 1.

Si $\delta =1$ on peut modifier la valeur de $\xi^n$ sur le compl\'ementaire de $A_n$ par un champ unitaire $\eta_x^n\in H_x$ quelconque, et obtenir ainsi une suite $\tilde \xi^n$ presque invariante et domin\'ee (par la fonction constante $\eta+1$), dont la norme ne tend pas vers 0, et est \`a distance uniform\'ement $<1$ (pour la norme $L^2$) de $\xi^n \in V$, o\`u $\norm{\xi^n}=1$. Ceci contredisant les hypoth\`eses, on a $\delta = 0$. 

\bigskip

En particulier pour tout $\eta \geqslant 1$ on obtient $\mu(A_n^\eta)\to 0$, o\`u
\[
A_n^\eta = \{ \frac 1 \eta \leqslant \norm {\xi^n} \leqslant \eta\}
\]
(quitte \`a extraire).  Par proc\'ed\'e diagonal, on peut donc trouver une suite extraite de $(\xi^n)$, encore not\'ee $(\xi^n)$, de sorte que  $\mu(A_n) \to 0$, o\`u
\[
A_n = \{ \frac 1 n \leqslant \norm {\xi^n} \leqslant n\}.
\]

Consid\'erons la fonction $f_n$ d\'efinie par 
\[
f_n(x) = 0\ \mathrm{si} \ \norm {\xi_x^n} \leqslant \frac 1 n \ \mathrm{et} \  f_n(x) = \norm {\xi_x^n} ^2\ \mathrm{sinon}.
\]
Ainsi $f_n \in L^1(X)$ et $\norm {f_n} _1 \to 1$. (Nous nous ramenons ici \`a  un argument de Connes-Feldman-Weiss  \cite[page 441]{CFW81}). Notons que 
\[
L^1E(f_n)=\int_X \sum_{(x,y)\in \R} \abs{f_n(y) - f_n(x)} \nu(x\to y) d\mu(x) \to 0.
\]
En effet
\begin{eqnarray}
L^1E(f_n) &\leqslant& \frac 2 n  + \int_\R (\norm{\pi(x,y)\xi^n_y - \xi^n_x})(\norm{\xi^n_y}+\norm{\xi_x^n}) d\h_\nu(x,y)\nonumber\\
&\leqslant&  \frac 2 n + 2\cdot \sqrt{E(\xi^n)}\norm{\xi^n}_2\to 0\nonumber
\end{eqnarray}
o\`u la deuxi\`eme in\'egalit\'e r\'esulte de l'in\'egalit\'e de Cauchy-Schwarz. D\'eduisons en l'existence de suites de F\o lner \'evanescentes dans $K$. Soit $\1_a$ la fonction caract\'eristique de $[a,\infty[ \subset \RI$.  Pour $t,t'\geqslant 0$ on a
\[
t = \int_0^\infty \1_a(t) da\quad \mathrm{et} \quad \abs {t-t'} = \int_0^\infty \abs{\1_a(t) - \1_a(t')}da.
\]
Soit $\varepsilon >0$ et $n$ suffisament grand pour que 
\[
L^1E(f_n) < \eta\varepsilon \norm{f_n}_1
\]
et $\mu\{f_n=0\}\geqslant 1-\varepsilon$, o\`u $\eta>0$ est une constante telle que $\nu(x\to y) \geqslant \eta$ sur $K$. Posons $f=f_n$. D'apr\`es le th\'eor\`eme de Fubini on a, en notant $f_a= \1_a\cdot f$,
\[
\int_0^\infty \int_\R \abs{f_a(y)-f_a(x)}\nu(x\to y)d\h(x,y)da<\eta\varepsilon \iint f_a d\h_\nu da.
\]
Consid\'erons donc $a>0$ tel que 
\[
\int_\R \abs{f_a(y)-f_a(x)}\nu(x\to y)d\h(x,y)<\eta\varepsilon \int_\R f_a d\h_\nu
\]
et notons $\Omega = \{ f_a = 1\}$. On a 
\begin{eqnarray}
\int_\R \abs{f_a(y)-f_a(x)}d\h_\nu(x,y)&\geqslant&\int_\Omega \sum_{y\sim x} \abs{f_a(y)-f_a(x)}\nu(x\to y) d\mu(x)\nonumber\\
&=&\int_X \left(\sum_{x\sim y}\chi_\Omega(x)\nu(y\to x)\right)\nonumber\abs{f_a(y)-1}d\mu(y)\nonumber\\
&=&\int_{\del \Omega} \left(\sum_{x\sim y}\chi_\Omega(x)\nu(y\to x)\right)\nonumber\abs{f_a(y)-1}d\mu(y)\nonumber\\
&\geqslant& \eta \int_{\del \Omega} \nonumber\abs{f_a(y)-1}d\mu(y)\nonumber\\
&=& \eta \mu(\del_K \Omega)\nonumber
\end{eqnarray}
Ainsi
\[
\mu(\del_K \Omega) < \varepsilon \mu(\Omega) \quad \mathrm{et} \quad 0<\mu(\Omega)\leqslant \varepsilon.
\]
Donc $K$ contient des suites de F\o lner \'evanescentes, ce qui contredit les hypoth\`eses.
\end{dem}

\bigskip

\begin{dfn}
Soit $\Q$ un espace singulier. On appelle \textnormal{marche al\'eatoire $\Q$-p\'eriodique} (sym\'etrique, born\'ee) la donn\'ee d'une d\'esingularisation discr\`ete $p:X\to \Q$ de $\Q$ et d'une marche al\'eatoire (sym\'etrique, born\'ee) $\nu$ sur les orbites de la relation $\R_p$. 

Plus pr\'ecis\'ement, la marche al\'eatoire $\Q$-p\'eriodique associ\'ee \`a cette donn\'ee est par d\'efinition la marche al\'eatoire $\tilde \nu$ sur l'ensemble d\'enombrable $\Q$-p\'eriodique $\R_p$ d\'efinie par
 \[
\nu((x,y)\to (x,z))=\nu(y\to z)
\]
(i.e. $\tilde \nu$ est l'extension \'equivariante de $\nu$ \`a $\R_p$). 
\end{dfn}

\bigskip

\textit{Remarque.} De m\^eme que dans \cite{Pichot04_I}, $\nu$ et $\tilde \nu$ correspondent aux deux points de vue transverse et longitudinal pour un m\^eme objet.

\bigskip
\bigskip

Passons \`a la d\'emonstration des th\'eor\`emes spectraux.

\bigskip

\begin{dfn}
Soit $\R$ une relation d'\'equivalence mesur\'ee sur un espace de probabilit\'e $(X,\mu)$, et $\nu$ une marche al\'eatoire sur cette relation. La \textnormal{diffusion simple} associ\'ee $\nu$ est l'op\'erateur de diffusion $D_\nu=D_{\nu,\1}$, agissant sur $L^2(X)$, associ\'e \`a la repr\'esentation triviale de $\R$. 
\end{dfn}

\bigskip

\begin{thm}
Soit $\Q$ un espace singulier ergodique de type fini. Si $\Q$ est de type $\II$, il poss\`ede un quotient moyennable si et seulement si les diffusions simples associ\'ees aux marches al\'eatoires $\Q$-p\'eriodiques sym\'etriques born\'ees n'ont pas de trou spectral.
\end{thm}

\bigskip

Ce th\'eor\`eme r\'esulte imm\'ediatement des r\'esultats de \cite{Pichot04_I} (cf. th. 14) et du th\'eor\`eme suivant, qui est \`a rapprocher des r\'esultats obtenus par K. Schmidt dans \cite{Schmidt81} pour des actions $\IIi$ de groupes d\'enombrables (cf. prop. 2.3). Nous renvoyons \'egalement ici \`a \cite{HjorthKechris03}.

\bigskip

\begin{thm}
Soit $K$ un graphage sym\'etrique born\'e d'une relation d'\'equivalence ergodique sur un espace de probabilit\'e $(X,\mu)$. On note $\nu=\nu_K$ la marche al\'eatoire associ\'ee \`a $K$. Alors $K$ contient des suites de F\o lner \'evanescentes (cf. \cite{Pichot04_I}) si et seulement si $\kappa_\nu=1$.
\end{thm}

\begin{dem}
Montrons d'abord que si $K$ contient des suites de F\o lner \'evanescentes, alors $\kappa_\nu=1$. Soit $A_n$ une suite de F\o lner dans $K$ et $\overline A_n = A_n \amalg \del_K A_n$. Consid\'erons l'espace $H$ des fonctions orthogonales aux constantes.  Alors $f_n =  \chi_{\overline A_n}-\mu(\overline A_n)$ de $\chi_{\overline A_n}$ sur $H$ est une suite non triviale de points presque fixes pour $D_\nu$, cf \cite{Schmidt81,Rosenblatt81,HjorthKechris03} (cette observation est d\^ue \`a K. Schmidt). Explicitement,
\[
\langle D_K f_n \mid f_n \rangle=\int_X\sum_y\chi_{\overline A_n}(x)\chi_{\overline A_n}(y)\nu(x\to y) d\mu(x)-\mu(\overline A_n)^2\geqslant \mu(A_n)-\mu(\overline A_n)^2
\]
et $\norm{f_n}^2=\mu(\overline A_n) - \mu(\overline A_n)^2 $. 

\bigskip

R\'eciproquement supposons que $K$  ne contienne pas de suites de F\o lner \'evanescentes, et montrons que la marche al\'eatoire $\nu$ associ\'ee \`a $K$ v\'erifie $\kappa_\nu<1$. V\'erifions les hypoth\`eses du lemme \ref{CFW} (\textit{iii}.). Soit $\xi_n$ une suite presque invariante domin\'ee et $\zeta^n$ le champ d\'efini par
\[
\zeta_x^n = \xi_x^n -\int_X \xi_x^nd\mu.
\]
Alors $\zeta^n \to 0$ presque s\^urement quitte \`a extraire (cf. lemme \ref{triv}), et le th\'eor\`eme en r\'esulte. 
\end{dem}

\bigskip

Nous pouvons maintenant achever la d\'emonstration du lemme \ref{triv}.

\bigskip

\begin{lem}\label{trivinf}
Soit $\R$ une relation d'\'equivalence fortement ergodique de type fini sur un espace de probabilit\'e $(X,\mu)$. On suppose que $\mu$ est une mesure de probabilit\'e invariante et on consid\`ere une suite $(\xi^n : X \to H)_n$ presque invariante pour la repr\'esentation triviale de $\R$ sur un espace de Hilbert $H$, domin\'ee par une fonction $g\in L^1(X)$, et telle que
\[
\norm{\int_X \xi^nd\mu}\to_n 0.
\]
Il existe une suite extraite de  $(\xi^n)$ qui converge presque surement vers 0. 
\end{lem}

\begin{dem} Soit $\xi^n$ une suite presque invariante domin\'ee. $\R$ \'etant une relation d'\'equivalence de type $\IIi$ fortement ergodique de type fini, elle poss\`ede un graphage  $K$  sym\'etrique born\'e ne contenant pas de suites de F\o lner \'evanescentes (cf. \cite{Pichot04_I}).  La marche al\'eatoire $\nu=\nu_K$ associ\'ee \`a $K$ d\'efinit une diffusion $D=D_\nu$ agissant sur $L^2(X)$ et poss\'edant un trou spectral. Soit $H$ un espace de Hilbert et $D^H=D_\nu^H$ la diffusion associ\'ee \`a la marche al\'eatoire $\nu$ et la repr\'esentation triviale de $\R$ sur $X\times H$. Alors $D$ poss\`ede un trou spectral si et seulement si $D^H$ en poss\`ede un (comme on le voit  par exemple \`a l'aide des in\'egalit\'es de Poincar\'e). Ainsi le lemme \ref{CFW} \textit{ii.} montre qu'il existe (quitte \`a extraire) une suite $\zeta^n$ de champs invariants  tels que 
\[
\norm {\xi_x^n -  \zeta^n_x}_x \to 0
\]
pour presque tout $x\in X$. Les champs invariants sont essentiellement constants, notons $\zeta_n\in H$ la valeur essentielle de $\zeta^n$. Alors l'hypoth\`ese
\[
\norm{\int_X \xi^nd\mu}\to_n 0
\]
entra\^\i ne que $\norm{\zeta_n}_H \to 0$ et le lemme en r\'esulte.
\end{dem}

\bigskip
\bigskip

Concluons ce paragraphe par la preuve du th\'eor\`eme \ref{T-intr}. Rappelons que les diffusions hilbertiennes associ\'ees \`a une marche al\'eatoire $\nu$ sur les orbites d'une relation d'\'equivalence mesur\'ee sont les op\'erateurs $D_{\nu,\pi}$ d\'efinis au paragraphe pr\'ec\'edent (associ\'es \`a chaque repr\'esentation hilbertienne $\pi$ de cette relation d'\'equivalence).

\bigskip

\begin{thm} 
Soit $\Q$ un espace singulier ergodique. Alors $\Q$ poss\`ede la propri\'et\'e T de Kazhdan si et seulement s'il existe une marche al\'eatoire $\Q$-p\'eriodique sym\'etrique dont les diffusions hilbertiennes poss\`edent un trou spectral.
\end{thm}

\begin{dem}
Soit $\Q$ un espace singulier ergodique poss\'edant la propri\'et\'e T. Alors $\Q$ est de type fini et poss\'ede une d\'esingularisation de type $\IIi$ sur un espace de probabilit\'e $(X,\mu)$. Cette d\'esingularisation ayant elle-m\^eme la propri\'et\'e T, elle poss\`ede un graphage $K$ sym\'etrique u.l.b. ne contenant pas de suites de F\o lner \'evanescentes. Soit $\nu=\nu_K$ la marche al\'eatoire r\'eguli\`ere associ\'ee \`a $K$, qui est sym\'etrique born\'ee relativement \`a $\mu$. Soit $\pi$ une repr\'esentation de $\R$. Le th\'eor\`eme \ref{prox} montre  que pour toute suite presque invariante $\xi^n$ telle que $\norm{\xi^n_x}_x\leqslant g(x)$ pour une fonction $g\in L^2(X)$, il existe (\`a extraction pr\`es) une suite $\zeta_n$ de champs invariants tels que 
\[
\norm {\xi_x^n -  \zeta^n_x}_x \to 0
\]
pour presque tout $x\in X$, et d'apr\`es le lemme \ref{CFW} \textit{iii.}, on a $\kappa_\pi < 1$.

\bigskip

R\'eciproquement supposons qu'il existe une marche al\'eatoire sym\'etrique sur une d\'esingularisation $p : X \to \Q$ de $\Q$ dont les diffusions hilbertiennes poss\`edent un trou spectral. Soit $\pi$ une repr\'esentation de $\R_p$ contenant une suite $\xi^n$ presque invariante domin\'ee telle que $\norm{\xi^n}_2 =1$. Le lemme \ref{CFW} \textit{ii.} entra\^ine qu'il existe une suite $\zeta_n$ de champs invariants tels que 
\[
\norm {\xi^n -  \zeta^n}_2 \to 0.
\]
Comme $\norm{\xi^n}_2=1$, on obtient $\norm{\zeta^n}_2 \neq 0$ pour $n$ grand, donc $\pi$ contient des champs invariants non triviaux. Donc $\R$ a la propri\'et\'e T (et $\Q$ \'egalement).
\end{dem}

\bigskip

\textit{Remarques.} 1 - Observons que si $\Q$ poss\`ede la propri\'et\'e T, on peut choisir une marche al\'eatoire $\Q$-p\'eriodique satisfaisant \`a la conclusion du th\'eor\`eme qui soit sym\'etrique born\'ee.

2 -  Soient $\R$ une relation d'\'equivalence mesur\'ee et $\nu$ une marche al\'eatoire sur ses orbites.  On d\'efinit  
\[
\kappa(\R, \nu) = \sup_\pi \kappa_{\nu,\pi},
\]
et de m\^eme $c_n(\R,\nu)$. Il r\'esulte ainsi du th\'eor\`eme ci-dessus que  s'il existe une marche al\'eatoire $\nu$ sym\'etrique born\'ee (relativement \`a une mesure de probabilit\'e $\mu$) telle que $c_2(\R,\nu) <2$, alors $\R$ a la propri\'et\'e T de Kazhdan (l'hypoth\`ese d'ergodicit\'e n'ayant pas \'et\'e utilis\'ee pour cette implication).  Ce r\'esultat sera utilis\'e sous cette forme au \textsection\ref{Garland} qui suit. Notons \'egalement que, de m\^eme que dans le cas des groupes, les constantes $\kappa_{\nu,\pi}$ et  $c_n(\R,\pi)$ sont uniformes en $\pi$ lorsque $\R$ poss\'ede la propri\'et\'e T.

\vspace{1cm}
\section{Le crit\`ere $\lambda_1 > 1/2$}\label{Garland}

Soit $\Q$ un espace singulier. Dans ce paragraphe nous donnons un crit\`ere, portant sur la structure locale des complexes simpliciaux $\Q$-p\'eriodiques de dimension 2, qui entra\^\i ne s'il est satisfait que l'espace singulier $\Q$ est un espace de Kazhdan. L'origine de ce crit\`ere se situe dans les travaux de Garland sur la cohomologie de groupes d'automorphismes d'immeubles de Bruhat-Tits cocompacts \cite{Garland73,Borel74}. Nous reprenons ici la d\'emonstration de ce r\'esultat donn\'ee par Gromov \cite{Gromov03}. 

\bigskip

Suivant Garland, cette d\'emonstration s'effectue en deux \'etapes,

1) l'\'etude spectrale de la g\'eom\'etrie locale (des \og links\fg) du complexe simplicial.

2) \og l'int\'egration\fg\ des donn\'ees spectrales locales en une in\'egalit\'e de Poincar\'e globale entra\^\i nant la propri\'et\'e T (\'etape dite de \og g\'eom\'etrie int\'egrale\fg\ dans \cite{Gromov03}).

\bigskip

Le concept d'\'energie introduit par Gromov \cite{Gromov03} dans ce contexte est essentiel pour l'\'etape 2 : il est lin\'eaire relativement \`a l'op\'eration d'addition (moyennisation) des diffusions (i.e. des marches al\'eatoires) et, par suite, les in\'egalit\'es de Poincar\'e s'int\`egrent. La premi\`ere \'etape  est classique. Une source importante de poly\`edre satisfaisant au crit\`ere local est donn\'ee par certains immeubles de Bruhat-Tits (de rang 2). Les propri\'et\'es spectrales de leurs links (des immeubles sph\'eriques) ont \'et\'e \'etudi\'ees par Feit-Higman \cite{FeitHigman64}.

\bigskip
\bigskip

Passons maintenant \`a la d\'emonstration du th\'eor\`eme annonc\'e (dont nous rappelons l'\'enonc\'e ci-dessous). La proc\'edure d\'ecrite ci-dessus n\'ecessite d'\'etudier le comportement de la marche al\'eatoire, sur le complexe simplicial, et sur ses links ; nous commen\c cons par les links.

\bigskip

\subsection{La condition $\lambda_1 > 1/2$.} Soit $L$ un graphe fini non orient\'e. On note $\tau(y)$ la valence d'un sommet $y$ dans $L$. La marche al\'eatoire uniforme sur $L$ est donn\'ee par 
\[
\nu (y\to z) = 1/\tau(y)
\]
si $(y,z) \in L^{(1)}$ est une ar\^ete de $L$ et $\nu(y\to z) =0$ sinon. Elle est sym\'etrique relativement \`a la mesure stationnaire $\mu(y)= \tau(y)/2\tau$, o\`u $\tau$ est le nombre d'ar\^etes de $L$. En d'autres termes on a
\[
\mu(y)\nu(y\to z) = \mu(z)\nu(z\to y)
\]
pour tout $y,z\in X$.

\bigskip

On note $D_L : \ell^2(L^{(0)},\mu) \to  \ell^2(L^{(0)},\mu)$ la diffusion associ\'ee, agissant sur les fonctions complexes d\'efinies sur les sommets de $L$, selon la formule
\[
D_Lf(y) = \frac 1 {\tau(y)}\sum_{z\sim y} f(z).
\]
Le produit scalaire sur $\ell^2(L^{(0)},\mu)$ est donn\'e par 
\[
\langle f\mid g\rangle = \sum_{y\in L} f(y)\overline {g(y)} \mu(y).
\]
$D_L$ est un op\'erateur hermitien et le laplacien associ\'e $\Delta_L = \id - D_L$ est positif, de noyau les fonctions constantes (on suppose $L$ connexe).\\

Le \textit{spectre de $L$} est par d\'efinition l'ensemble des valeurs propres $\Delta_L$. Le \textit{trou spectral} de $L$ est la plus petite valeur propre non nulle et se note $\lambda_1(L)$.\\

Soit $H$ un espace de Hilbert. On \'etend $D_L$ en un op\'erateur born\'e encore not\'e $D_L$ agissant sur $\ell^2(L^{(0)},H,\mu)$. Il est hermitien pour la norme hilbertienne et sans point fixe sur l'orthogonal des constantes ; sa plus grande valeur propre est encore $\kappa= 1 - \lambda_1(L)$. 

\bigskip

Soit $\nu_\infty(y\to z) =  \mu(z)$ la marche al\'eatoire stationnaire associ\'ee \`a $\nu$. L'in\'egalit\'e de Dirichlet-Poincar\'e correspondante s'\'ecrit
\[
E_\infty(\xi) \leqslant c_\infty E(\xi)
\]
o\`u
\[
c_\infty = 1 + \kappa + \kappa^2 + \ldots = \frac 1 {\lambda_1(L)}.
\]
En d'autres termes on a,
\[
\frac 1 {2\tau} \sum_{y,z\in L^{(0)}} \norm{\xi_z -\xi_y}^2 \tau(y)\tau(z) \leqslant \frac 1 {\lambda_1} \sum_{(y,z)\in L^{(1)}}\norm{\xi_z -\xi_y}^2,
\]
pour toute fonction $\xi \in \ell^2(L^{(0)},H,\mu)$.\\

\subsection{G\'eom\'etrie locale et int\'egrale.}

Soit $\Q$ un espace singulier de type fini. On consid\'ere un complexe simplicial $\Q$-p\'eriodique u.l.f. $\tilde \Sigma$ de dimension 2 et on note $\Sigma\to \Q$ la d\'esingularisation simpliciale associ\'ee \`a $\tilde \Sigma$.  L'ensemble $X$ des sommets de $\Sigma$ est un espace bor\'elien standard muni d'une relation d'\'equivalence $\R$ \`a classes d\'enombrables, qui d\'etermine une d\'esingularisation discr\`ete de $\Q$ ; \'etant donn\'ee une mesure de probabilit\'e quasi-invariante $\mu$ sur $X$, on note $\delta$ le cocycle de Radon-Nikodym associ\'e \`a $\mu$.

\bigskip

\begin{dfn} Le \textnormal{link} en un sommet $s$ de $\Sigma$ est le graphe fini dont les sommets sont les ar\^etes de $\Sigma$ issues de $s$ et les ar\^etes les triangles de $\Sigma$ issues de $s$.
\end{dfn}

\bigskip

\begin{thm} Soit $\Q$ un espace singulier de type fini, $\tilde \Sigma$ un complexe simplicial $\Q$-p\'eriodique u.l.f.  de dimension 2, et $\Sigma$ la d\'esingularisation simpliciale associ\'ee. On fixe une mesure quasi-invariante $\mu$ sur la transversale des sommets $X=\Sigma^{(0)}\subset \Sigma$ et on consid\`ere un nombre r\'eel $\delta_\mu \geqslant 1$ tel que 
\[
\delta_\mu^{-1} \leqslant \delta(y,z) \leqslant \delta_\mu
\]
pour presque toute ar\^ete $(y,z)$ de $\Sigma^{(1)}$. On suppose que presque tout link $L$ de $\Sigma$ est connexe et v\'erifie $\lambda_1(L) \geqslant \lambda$, o\`u $\lambda$ est un nombre r\'eel v\'erifiant
\[
\lambda > \delta_\mu^3/2.
\]
Alors $\Q$ poss\`ede la propri\'et\'e T de Kazhdan.
\end{thm}

\begin{dem} On note $K \subset \R$ le graphage de $\R$ associ\'e au 1-squelette $\Sigma^{(1)}$ de $\Sigma$.  Si $x\in X$ est un sommet de $\Sigma$, on note $L_x$ son link dans $\Sigma$, et $\lambda_1(x)$ la plus petite valeur propre non nulle de $L_x$. Enfin, on note $\tau(y,z)$ le nombre de triangles de $\Sigma$ contenant l'ar\^ete $(y,z) \in \Sigma^{(1)}$, et  $\tau(x)$ le nombre de triangles attach\'es au sommet $x\in X$ dans $\Sigma$.

\bigskip

Supposons $\lambda_1(x) \geqslant \lambda$ et consid\'erons une repr\'esentation $\pi$ de $\R$. Soit $\xi : X \to H$ un champ de carr\'e int\'egrable.  
 Par hypoth\`ese on a, en notant $\overline \xi$ l'extension \'equivariante de $\xi$ \`a $\R$ (d\'efinie par $\overline \xi_x(y) = \pi(x,y)\xi_y$ pour $x\sim_\R y$),
\[
E_{\underline \nu_x}(\overline \xi_x) \leqslant \frac 1 \lambda E_{\nu_x}(\overline \xi_x)
\]
o\`u $\nu_x (y\to z) = 1/\tau(x,y)$ est la marche al\'eatoire uniforme sur le graphe fini $L_x$, sym\'etrique relativement \`a la mesure stationnaire $\mu_x(y) = \tau(x,y)/2\tau(x)$, et $\underline \nu_x(y\to z) =\tau(x,z)/2 \tau(x)$ est la marche al\'eatoire stationnaire correspondante. Explicitement, cette in\'egalit\'e s'\'ecrit 
\[
\frac 1 {2\tau(x)} \sum_{y,z \in L_x^{(0)}} \norm{d\xi(y,z)}^2\tau(x,y)\tau(x,z) \leqslant \frac 1 \lambda \sum_{(y, z)\in L_x^{(1)}}\norm{d\xi(y,z)}^2.
\]

\bigskip

\noindent Par ailleurs on a
\begin{eqnarray}
\int_X  \sum_{(y, z)\in L_x^{(1)}}\norm{d\xi(y,z)}^2 d\mu(x)&=& \int_X \sum_{(y,z)\in K} \sum_{(x,y,z)\in \Sigma }\norm{d\xi(y,z)}^2 \delta(x,y)d\mu(y)\nonumber \\ 
&=&\int_X \sum_{(y,z)\in K} \norm{d\xi(y,z)}^2 \tau_\delta(y,z)d\mu(y)\nonumber
\end{eqnarray}
o\`u
\[
\tau_\delta(y,z) = \sum_{(x,y,z)\in \Sigma} \delta(x,y)
\]
(la somme porte sur les triangles contenant $(y,z)$), et
\begin{eqnarray}
&&\int_X\frac 1 {2\tau(x)} \sum_{y,z \in L_x^{(0)}} \norm{d\xi(y,z)}^2\tau(x,y)\tau(x,z)d\mu(x)\nonumber\\
&&\qquad\qquad=\int_X \sum_{z\sim y}\norm{d\xi(y,z)}^2\sum_{L_x^{(0)} \ni y,z} \frac 1 {2\tau(x)}\tau(x,y)\tau(x,z)\delta(x,y)d\mu(y)\nonumber\\
&&\qquad\qquad = \int_X \sum_{z\sim y} \norm{d\xi(y,z)}^2 \underline\tau_\delta(y,z)d\mu(y)\nonumber
\end{eqnarray}
o\`u $\underline \tau_\delta(y,z) = \sum_{L_x^{(0)} \ni y,z} \frac 1 {2\tau(x)}\tau(x,y)\tau(x,z)\delta(x,y)$. Posons
\[
\tau_\delta(y)=\frac 1 2 \sum_{(x,y)\in K} \delta(x,y)\tau(x,y).
\]
On a
\[
\sum_{z\sim y} \tau_\delta(y,z) = \sum_{z\sim y} \underline \tau_\delta(y,z)= 2\tau_\delta(y),
\]
de sorte que, en d\'efinissant
\[
\nu(y\to z) =\frac{\tau_\delta(y,z)}{2\tau_\delta(y)}
\]
et
\[
\underline \nu(y\to z) =\frac{\underline \tau_\delta(y,z)}{2\tau_\delta(y)}
\]
on obtient une in\'egalit\'e de Poincar\'e,
\[
\int_X \sum_{(y,z)\in \R}  \norm{d\xi(y,z)}^2 \underline\nu(y\to z) d\tilde\mu(y) \leqslant \frac 1 \lambda \int_X\sum_{(y,z)\in\R}\norm{d\xi(y,z)}^2 \nu(y\to z) d\tilde\mu(y)
\]
o\`u $d\tilde \mu(y) = 2\tau_\delta(y) d\mu(y)$. Remarquons que $\tilde \mu$ est sym\'etrique relativement \`a $\nu$ et $\underline \nu$.\\

Il est facile de voir que si $\mu$ est invariante, alors $\underline \nu =\nu^2$ est la marche al\'eatoire en deux pas associ\'ee \`a $\nu$ (ainsi $\R$ poss\`ede la propri\'et\'e T si $\lambda > 1/2$, comparer \`a \cite{Zuk96,Pichot03}). Sinon on a
\[
\nu^2(y\to z) = \sum_{L_x^{(0)}\ni y,z} \frac{\tau_\delta(y,x)}{2\tau_\delta(y)}\frac{\tau_\delta(x,z)}{2\tau_\delta(x)}.
\]
Par hypoth\`ese $\tau_\delta(x,y) \leqslant\delta_\mu \tau(x,y)$, $\tau_\delta(x,z) \leqslant\delta_\mu \tau(x,z)$ et $\tau_\delta(x)\geqslant \tau(x)/\delta_\mu$ ; donc
\[
\nu^2(y\to z) \leqslant \delta_\mu^3 \underline \nu(y\to z)
\]
et on obtient l'in\'egalit\'e de Poincar\'e
\[
E_{\nu^2}(\xi) \leqslant \frac {\delta_\mu^3} \lambda E_\nu(\xi).
\]
Donc $\R$ (ainsi que $\Q$) a la propri\'et\'e T si $\lambda > \delta_\mu^3/2$.
\end{dem}

\bigskip
\bigskip

\begin{cor}
Le 1-squelette d'un complexe simplicial $\Q$-p\'eriodique satisfaisant aux hypoth\`eses du th\'eor\`eme ne contient pas de suites de F\o lner \'evanescentes.
\end{cor}

\bigskip

En effet la diffusion simple associ\'ee \`a la marche al\'eatoire sur le 1-squelette de ce complexe poss\`ede un trou spectral.

\bigskip


\subsection{L'espace des immeubles $\tilde A_2$.} Il existe une infinit\'e non d\'enombrable d'immeubles $\tilde A_2$ d'\'epaisseur fix\'ee. Au cours d'une discussion avec \'Etienne Ghys, nous avons constat\'e qu'il \'etait possible de construire une lamination compacte dont l'espace des feuilles coincide avec l'ensemble de ces immeubles. L'existence de cette lamination fut l'une des motivations de cet article. Le link d'un immeuble de type $\tilde A_2$ v\'erifie le crit\`ere $\lambda_1 >1/2$, mais l'existence de mesures invariantes sur cette lamination n'est pas claire \textit{a priori}.  \'Etant donn\'e que nous entreprenons une \'etude d\'etaill\'ee de cet espace dans  \cite{BarrePichot04_I,BarrePichot04_II}, en collaboration avec Sylvain Barr\'e, il est inutile que nous en disions davantage ici.

\bigskip
\bigskip
\bigskip
\noindent  Unit\'e de Math\'ematiques Pures et Appliqu\'ees,\\
Unit\'e Mixte de Recherche CNRS 5669\\
46, all\'ee d'Italie,\\
69364 Lyon Cedex 07, France.\\
mpichot@umpa.ens-lyon.fr, pichot@ms.u-tokyo.ac.jp.

\end{document}